\documentclass[11pt]{article}

\usepackage{amsfonts, amssymb, amsmath}
\begin{document}
\newtheorem{thm}{Theorem}[section]
\newtheorem{conj}{Conjecture}[section]
\newtheorem{lm}[thm]{Lemma}
\newtheorem{prop}[thm]{Proposition}
\newtheorem{cor}[thm]{Corollary}
\newcommand{\sign}{\mathrm{\mathop{sign}}}
\newcommand{\w}{\widehat}
\newcommand{\D}{\mathcal{D}}
\newcommand{\QED}{\hfill \blacksquare}
\newcommand{\proof}{\noindent {\bf Proof.}\ }
\renewcommand{\baselinestretch}{1.1}
\newtheorem{rem}{Remark}[section]

\title{\bf  Eigenvalues of the Derangement Graph}
\author{
Cheng Yeaw Ku
\thanks{ Department of Mathematics, Caltech, Pasadena, CA 91125,
USA. E-mail: cyk@caltech.edu.} \and David B. Wales \thanks{
Department of Mathematics, Caltech, Pasadena, CA 91125, USA. E-mail:
dbw@caltech.edu.} }

\maketitle

\begin{abstract}
We consider the Cayley graph on the symmetric group $S_{n}$
generated by derangements. It is well known that the eigenvalues of
this graph are indexed by partitions of $n$. We investigate how
these eigenvalues are determined by the shape of their corresponding
partitions. In particular, we show that the sign of an eigenvalue
is the parity of the number of cells below the first row of the
corresponding Ferrers diagram. We also provide some lower and upper bounds
for the absolute values of these eigenvalues.

\bigskip\noindent
{\sc keywords:}
derangement, Cayley graph, eigenvalue, symmetric group


\end{abstract}

\section
{Introduction}\label{intro}

Let $G$ be a finite group and let $S$ be a nonempty subset of $G$
satisfying the condition that $s \in S \Longrightarrow s^{-1} \in S$
and $1 \not \in S$.  The {\em Cayley graph} $\Gamma(G,S)$ has the
elements of $G$ as its vertices and two vertices $u, v \in G$ are
joined by an edge if and only if $uv^{-1} \in S$.  In this 
paper, we shall be interested in the graph
$\Gamma_{n}$ which is $\Gamma(S_{n},\D_{n})$ where $S_{n}$ is the symmetric
group of permutations of the integers $1, \ldots, n$, denoted $[n]$, and $\D_n$ is 
the set of derangements of $[n]$ which are the permutations in $S_n$ which 
fix no point, i. e. for 
which $g(x) \not = x$ for all $ x\in [n]$. The graph $\Gamma_{n}$ is called the {\em derangement
graph} on $[n]$. Clearly, $\Gamma_{n}$ is vertex-transitive and so
it is $D_{n}$-regular, where $D_{n}=|\D_{n}|$. By a standard result
in graph theory, $D_{n}$ is the largest eigenvalue of $\Gamma_{n}$ see \cite{Gods}.

For a graph $\Gamma$, let $\alpha(\Gamma)$ denote the {\em
independence number} of $\Gamma$, i.e. the cardinality of an
independent set of maximum size of $\Gamma$.  For any $k$-regular
graph $\Gamma$ with $N$ vertices, the independence number satisfies
the Delsarte-Hoffman bound:
\[ \alpha(\Gamma) \le \frac{-\mu N}{k-\mu}, \]
where $\mu$ is the smallest eigenvalue of $\Gamma$, see \cite{Hoff}.  In particular,
this implies that $\mu$ is negative .

Without applying the Delsarte--Hoffman bound, Deza and Frankl
\cite{DF} first proved that $\alpha(\Gamma_{n})=(n-1)!$ by purely
combinatorial means. More recently, the structure of maximum-size
independent sets of $\Gamma_{n}$ has been determined by several
authors (\cite{CK}, \cite{GM}, \cite{LM}, \cite{WZ}) using different
methods, namely such a set must be a coset of the stabilizer of a
point. Ku and Wong \cite{KW} conjectured that $-\frac{D_{n}}{n-1}$
is the smallest eigenvalue of $\Gamma_{n}$ which would give
equality to the Delsarte--Hoffman bound. This has been proved by
Renteln \cite{Rent}.  It is immediate from this that $\alpha(\Gamma)=(n-1)!$
as the stabilizer of a point is an independent set.

The main theme of this work is to describe in more detail the
properties of the eigenvalues of the derangement graph. Recall that
a Cayley graph $\Gamma(G, S)$ is {\em normal} if $S$ is closed under
conjugation. Its spectra is described in the following lemma.  See
for example Lubotsky \cite[Theorem~$8.2.1$]{Lub}.  
\begin{lm}
The eigenvalues of a normal Cayley graph $\Gamma(G,S)$ are integers
given by
\begin{eqnarray}
\eta_{\chi} & = & \frac{1}{\chi(1)} \sum_{s \in S} \chi(s),
\nonumber
\end{eqnarray}
where $\chi$ ranges over all the irreducible characters of $G$.
Moreover, the multiplicity of $\eta_{\chi}$ is $\chi(1)^{2}$.
\end{lm}

Since $\D_{n}$ is closed under conjugation, $\Gamma_{n}$ is normal.
It is well known that both the conjugacy classes of $S_{n}$ and the
irreducible characters of $S_{n}$ are indexed by partitions
$\lambda$ of $n$ \cite{James}. Recall that a {\em partition} $\lambda$ of $n$,
denoted by $\lambda \vdash n$, is a weakly decreasing sequence
$(\lambda_{1}, \ldots, \lambda_{r})$ with $\lambda_{r} \ge 1$ such
that $\sum_{i=1}^{r} \lambda_{i} = n$. Its {\em size} is
$|\lambda|$, its {\em length} is $r$ and each $\lambda_{i}$ is the
$i$-{\em th part} of the partition. We also adopt the notation
$(\mu_{1}^{a_{1}}, \mu_{2}^{a_{2}}, \ldots, \mu_{s}^{a_{s}}) \vdash
n$ where $\mu_{i}$ are the distinct non zero parts which occur with
multiplicity $a_{i}$. For example,
\[ (5,4,4,3,3,3,1) \longleftrightarrow (5, 4^{2}, 3^{3},1). \]

Let $\lambda = (\lambda_{1}, \ldots, \lambda_{l})$ and
$\mu=(\mu_{1}, \ldots, \mu_{m})$ be partitions of $n$. Then $\lambda
> \mu$ in {\em lexicographic order} if, for some index $i$,
\[ \lambda_{j}=\mu_{j} ~~\textnormal{for } j<i ~~\textnormal{and }
\lambda_{i} > \mu_{i}. \] Note that for $\lambda_{1} \ge
\frac{n}{2}$, the partition $(\lambda_{1}, n-\lambda_{1})$ is the
largest partition in lexicographic order among the partitions with
the same first part, $\lambda_1$.

Based on the remarks above, we write $\eta_{\lambda}$ to denote the
eigenvalue $\eta_{\chi_{\lambda}}$ of $\Gamma_{n}$, where
$\chi_{\lambda}$ is the irreducible  character  indexed by the
partition $\lambda \vdash n$. We shall investigate how these
eigenvalues are determined by the shape of their corresponding
partitions. We used \textsf{GAP} \cite{GAP} to list the eigenvalues
for many values of $n$. Some values for small $n$ are tabulated in
Section~\ref{values}. A glance at the table shows many striking
properties about the values.  We prove the following main results and
offer a conjecture.

\begin{thm}[The Alternating Sign Property (ASP)]\label{ASP}
For any partition $\lambda=(\lambda_{1}, \ldots, \lambda_{r}) \vdash
n$,
\begin{eqnarray}
\sign(\eta_{\lambda}) & = & (-1)^{n-\lambda_{1}} \nonumber \\
& = & (-1)^{\textnormal{\# cells under the first row of }\lambda}
\label{1.1}
\end{eqnarray}
where $\sign(\eta_{\lambda})$ is $1$ if $\eta_{\lambda}$ is positive
or $-1$ if $\eta_{\lambda}$ is negative.
\end{thm}

\begin{thm}\label{Main2}
Let $\lambda = (\lambda_{1}, \ldots, \lambda_{r}) \vdash n$.

\noindent \textnormal{(i)} If $\lambda_{1} \ge \lfloor \frac{n}{2}
\rfloor$ then
\[ |\eta_{(\lambda_{1}, 1^{n-\lambda_{1}})}| \le |\eta_{\lambda}| \le |\eta_{(\lambda_{1}+1,
1^{n-\lambda_{1}-1})}|.
\]

\noindent \textnormal{(ii)} If $\lambda_{1} < \lfloor \frac{n}{2}
\rfloor$ then
\[  |\eta_{\lambda}| \leq  |\eta_{\big(\lfloor \frac{n}{2} \rfloor+1, 1^{n-\lfloor
\frac{n}{2}\rfloor-1}\big)}|, \] with strict inequality if $n\geq 6$.
Moreover, $|\eta_{(\lambda_{1}, 1^{n-\lambda_{1}})}| =
D_{\lambda_{1}}+(n-\lambda_{1})D_{\lambda_{1}-1}$ for any $\lambda_1$.  
Here as usual $\lfloor \frac{n}{2}\rfloor$ is the greatest integer less than or equal 
to $\frac{n}{2}$. 
\end{thm}

\begin{thm}\label{Main3}
Let $\lambda_{1}$ be $n-1$, $n-2$, $n-3$ or $n-4$ with $n$ being at least $2$, $4$, $6$, and 
$8$ respectively. Then the absolute values of the
eigenvalues which correspond to the partitions of $n$ with
$\lambda_{1}$ as their first part decrease in lexicographic order,
i.e. $|\eta_{\mu}| \le |\eta_{\lambda}|$ if and only if $\mu <
\lambda$ in lexicographic order. The decrease is strict for $n\ge 7$.  
\end{thm}

\begin{rem}  \rm Because of Theorem~\ref{ASP} there are no eigenvaulues~$0$ for $n\geq 2$ and 
so the adjacency matrix is nonsingular.  Note, however, that there is always 
an eigenvalue $(-1)^n$ from the partition $(2,1^{n-2})$ by Lemma~\ref{hook}.  
\end{rem}

\begin{rem} \rm Theorem \ref{Main2} implies that for $\lambda_{1} \ge\lfloor \frac{n}{2}
\rfloor$, $|\eta_{\lambda}| \le |\eta_{\mu}|$ whenever $\lambda_1 <
\mu_{1}$. In fact, we shall prove that the upper bound in part
(i) of Theorem \ref{Main2} is strict in most cases, namely when
$\lambda \not = (\frac{n}{2}, \frac{n}{2})$ or $\lambda \not =
(\frac{n-1}{2}, c, d)$ for any $c \le \frac{n-1}{2}$, $d \ge 1$ (see
Proposition~\ref{sqeezebelow-2}). Our approach does not
yield a good upper bound for $|\eta_{\lambda}|$ when $\lambda_{1}$
is small, i.e. $\lambda_{1} < \lfloor \frac{n}{2} \rfloor$.  The
condition $\lfloor \frac{n}{2} \rfloor \le \lambda_{1}$ for the
upper bound in part (i) of Theorem \ref{Main2} cannot be weakened
for otherwise there are many counterexamples, e.g. for $n=9$, we
have $|\eta_{(3^3)}|=32 > 19 = |\eta_{(4,1^5)}|$.
Nevertheless, part (ii) of Theorem \ref{Main2} implies that
$|\eta_{(3^3)}| < |\eta_{(5,1^4)}|=80$.  For $\lambda_{1} <
\frac{n}{2}$, it is generally not true that $|\eta_{\lambda_{*}}| <
|\eta_{(\lambda_{1}+1, 1^{n-\lambda_{1}-1})}|$ where $\lambda^{*}$
is the largest partition in lexicographic order among the partitions
with $\lambda_{1}$ as their first part.
\end{rem}

\begin{rem}\label{lexicographic} \rm Notice that Theorem~\ref{Main3} is a strengthening of Theorem~\ref{Main2} for partitions with
$\lambda_{1} \ge n-4$ and $n\geq 8$.
\end{rem}

\begin{rem} \rm It seems at first sight from our
computations that the absolute values of eigenvalues should decrease
in lexicographic order among the partitions with the same first part
$\lambda_{1}$ for all $\lambda_{1}$. However, this is not true in
general for both large and small values of $\lambda_{1}$ with respect to
$n$.  For example when
$n=15$, $|\eta_{(7,4,1^4)}|=5558<5566=|\eta_{(7,3^2,2)}|$ but in lexicographic order
$(7,4,1)>(7,3,2^2)$.  In fact, the smallest $n$ for which this occurs is $n=11$ with $|\eta_{(4,3,1^4)}|=37<38=|\eta_{(4,2^3,1)}|$.  Also, when
$n=17$ we have $|\eta_{(9,5,1^3)}| = 347104 < 349624 =
|\eta_{(9,4^2)}|$.  Notice here $\lambda_1>\frac{n}{2}$.  These values have been computed 
in GAP but are not in the tables in Section~\ref{values}.  
\end{rem}

In view of Theorem \ref{Main2}, Theorem \ref{Main3}, and values for small $n$ we make the
following conjecture:

\begin{conj}\label{Conjecture1}
Suppose $\lambda^{*}  \vdash n$ is the largest partition in
lexicographic order among all the partitions with $\lambda_{1}$ as
their first part. Then, for every $\lambda = (\lambda_{1}, \ldots,
\lambda_{s}) \vdash n$,
\[|\eta_{(\lambda_{1}, 1^{n-\lambda_{1}})}| \le  |\eta_{\lambda}| \le |\eta_{\lambda^{*}}|. \]
\end{conj}

Notice this follows for $n\geq 8$ and $\lambda_1\ge n-4 $ by Theorem~\ref{Main3}.  For $\lambda_1\ge \lfloor \frac{n}{2}\rfloor$
the lower bound holds by Theorem~\ref{Main2}.  The upper bound obtained in Theorem~\ref{Main2} is a different
upper bound.

Our main results can be regarded as a strengthening of the following
result of Renteln:

\begin{thm}[Renteln, \cite{Rent}] \label{Renteln_1}
The smallest eigenvalue of $\Gamma_{n}$ is
\[ \eta_{(n-1,1)} = -\frac{D_{n}}{n-1} = -(D_{n-1}+D_{n-2}). \]
Moreover, if $\lambda \not = (n), (n-1,1)$ then
$|\eta_{\lambda}|<|\eta_{(n-1,1)}|<|\eta_{(n)}|$.
\end{thm}
Notice Theorem~\ref{Main3} is an extension of this result.

Our methods rely heavily on the following remarkable recurrence
formula for the eigenvalues of $\Gamma_{n}$ proved by Renteln
\cite{Rent}. To describe this result, we require some terminology.
To the Ferrers diagram of a partition $\lambda$, we assign
$xy$-coordinates to each of its boxes by defining the
upper-left-most box to be $(1,1)$, with the $x$ axis increasing to
the right and the $y$ axis increasing downwards. Then the {\em hook}
of $\lambda$ is the union of the boxes $(x',1)$ and $(1, y')$ of the
Ferrers diagram of $\lambda$, where $x' \ge 1$, $y' \ge 1$. Let
$\w{h}_{\lambda}$ denote the hook of $\lambda$ and let $h_{\lambda}$
denote the size of $\w{h}_{\lambda}$. Similarly, let
$\w{c}_{\lambda}$ and $c_{\lambda}$ denote the first column of
$\lambda$ and the size of $\w{c}_{\lambda}$ respectively. Note that
$c_{\lambda}$ is equal to the number of rows of $\lambda$. When
$\lambda$ is clear from the context, we replace $\w{h}_{\lambda}$,
$h_{\lambda}$, $\w{c}_{\lambda}$ and $c_{\lambda}$ by $\w{h}$, $h$,
$\w{c}$ and $c$ respectively. Let $\lambda-\w{h} \vdash n-h$ denote
the partition obtained from $\lambda$ by removing its hook. Also,
let $\lambda-\w{c}$ denote the partition obtained from $\lambda$ by
removing the first column of its Ferrers diagram, i.e.
$(\lambda_{1}, \ldots, \lambda_{r})-\w{c} = (\lambda_{1}-1, \ldots,
\lambda_{r}-1) \vdash n-r$.

\begin{thm}[Renteln, \cite{Rent}] \label{Renteln_2}
For any partition $\lambda=(\lambda_{1}, \ldots, \lambda_{r})$,
\begin{eqnarray}
\eta_{\lambda} & = &
(-1)^{h}(\eta_{\lambda-\w{h}}+(-1)^{\lambda_{1}}h\eta_{\lambda-\w{c}})
\label{mainrecur}
\end{eqnarray}
with initial condition $\eta_{\emptyset}=1$.
\end{thm}
Since the above recurrence will be used extensively throughout the
paper, we shall often refer to it as the {\em main recurrence}.

The rest of the paper is organized as follows. In Section~\ref{preliminary}, we give
some useful formulae for the eigenvalues which correspond to
partitions of simple shapes. Theorem \ref{ASP} is proved in Section~\ref{alternatingsignprop}. In
Sections~\ref{lowerboundeta} and \ref{upperboundeta} respectively, we provide a lower and an upper
bound for the eigenvalues when the first part of the corresponding
partition is large. In Section~\ref{fewpartslambda} we consider partitions with few
parts. We prove Theorem \ref{Main2} in Section~\ref{secondupperboundeta}.  Sections~\ref{casen-2},
\ref{casen-3}, and \ref{casen-4}
are devoted to a proof of Theorem \ref{Main3}.  As mentioned, some values
for small $n$ have been tabulated in Section~\ref{values}.


\section{Some preliminary results}\label{preliminary}

In this section, we collect some basic formulae for some special
types of partitions. Recall the following useful facts about the
derangement numbers:

\begin{lm}\label{recurrence1}
For $n \ge 1$,
\begin{itemize}
\item[$(1)$] $D_{n} = nD_{n-1} + (-1)^{n}$.
\item[$(2)$] $D_{n} = (n-1)(D_{n-1}+D_{n-2})$.
\item[$(3)$] The first eleven derangement numbers are $D_0=1$, $D_1=0$, $D_2=1$, $D_3=2$, $D_4=9$, $D_5=44$,
$D_6=265$, $D_7=1854$, $D_8=14833$, $D_9=133496$, and $D_{10}=1334961$.
\end{itemize}
\end{lm}
\proof See \cite{Stanley}, page $67$ for $(1)$ and $(2)$.  The values in $(3)$ are tabulated in the first entries
in the tables in section~\ref{values} as $\eta_{\lambda}$ for $\lambda=(n)$.
$\QED$

\begin{rem}\label{increasing}  \rm Notice it follows from Lemma~\ref{recurrence1} that the
values of $D_n$ are striclty increasing for $n\ge 1$.
\end{rem}

\begin{lm}\label{hook}
Let $\lambda=(\lambda_1, 1^{n-\lambda_1})$ be a hook. Then
\[ \eta_{(\lambda_1, 1^{n-\lambda_1})} = (-1)^{n}\left( 1 + (-1)^{\lambda_1}nD_{\lambda_1-1}\right)
  =(-1)^{n-\lambda_1}\big(D_{\lambda_1}+(n-\lambda_1)D_{\lambda_1-1}\big).
  \] In particular, $\eta_{(n-1,1)} = -\frac{D_{n}}{n-1} = -(D_{n-1}+D_{n-2})$.
\end{lm}
\proof The first equality is \cite[Lemma~$7.4$]{Rent}.

For the second use the first.
\begin{eqnarray*}
\eta_{(\lambda_1,1^{n-\lambda_1})} &=&  (-1)^{n}\left( 1 +
(-1)^{\lambda_1}nD_{\lambda_1-1}\right)  \cr
     &=& (-1)^{\lambda_1+(n-\lambda_1)}\left(1+(-1)^{\lambda_1}\big(\lambda_1D_{\lambda_1-1}+
                               (n-\lambda_1)D_{\lambda_1-1}\big)\right) \cr
         &=&(-1)^{n-\lambda_1}\left((-1)^{\lambda_1} +\lambda_1D_{\lambda_1-1}+(n-\lambda_1)D_{\lambda_1-1} \right) \cr
           &=& (-1)^{n-\lambda_1}\left((-1)^{\lambda_1}+D_{\lambda_1}-(-1)^{\lambda_1}+
                                     (n-\lambda_1)D_{\lambda_1-1}\right) \cr
           &=& (-1)^{n-\lambda_1}\big( D_{\lambda_1}+(n-\lambda_1)D_{\lambda_1-1}\big).
\end{eqnarray*} Applying this when $\lambda_{1}=n-1$ gives the last
statement. $\QED$

We define a partition $\lambda = (\lambda_{1}, 2,
1^{n-\lambda_{1}-2})$ to be a {\em near hook}.

\begin{lm}\label{nearhook}
Let $\lambda=(\lambda_1,2,1^{n-\lambda_1-2})$ be a near hook. Then
\begin{eqnarray}
\eta_{(\lambda_1,2,1^{n-\lambda_1-2})} & = & (-1)^{n+\lambda_1}(n-1)\frac{D_{\lambda_1}}{\lambda_1-1} \nonumber \\
 & = & (n-1)\left((-1)^{n-1} +
(-1)^{n+\lambda_1} \lambda_1 D_{\lambda_1-2}\right). \nonumber
\end{eqnarray}
\end{lm}

\proof Use the main recurrence and the properties above or \cite[Lemma~$8.3$]{Rent} . $\QED$

\begin{lm}\label{hookcasevalues}  The values of $|\eta_\lambda|$ for hooks are given by Lemma~\ref{hook}.  The absolute values are as follows.
\begin{eqnarray*}
|\eta_{(1^n)}|&=& n-1  \cr |\eta_{(2,1^{n-2})}| &=& 1  \cr |\eta_
{(3,1^{n-3})}| &=& n-1
\end{eqnarray*}
For $2 \leq \lambda_{1} \leq n-1$ we have
$|\eta_{(\lambda_{1},1^{n-\lambda_{1}})}|<|\eta_{(\lambda_{1}+1,1^{n-\lambda_{1}-1})}|.$
\end{lm}
\proof By Lemma~\ref{hook} we see
$|\eta_{(\lambda_{1},1^{n-\lambda_{1}})}|=D_{\lambda_{1}}+(n-\lambda_{1})D_{\lambda_{1}-1}$.
The values for $\lambda_{1} \leq 3$ are as given. We need only show
for $3\leq \lambda_{1} \leq n-1$ that
$|\eta_{(\lambda_{1},1^{n-\lambda_{1}})}|<
|\eta_{(\lambda_{1}+1,1^{n-\lambda_{1}-1})}|$.  This means we need
to show $D_{\lambda_{1}}+(n-\lambda_{1})D_{\lambda_{1}-1}<
D_{\lambda_{1}+1}+(n-\lambda_{1}-1)D_{\lambda_{1}}$ for these values
of $\lambda_{1}$. Using
$(n-\lambda_{1})D_{\lambda_{1}-1}<(n-\lambda_{1})D_{\lambda_{1}}$ we
need only show $0<D_{\lambda_{1}+1}-2D_{\lambda_{1}}$ which is true
using $D_{\lambda_{1}+1}=(\lambda_{1}+1)D_{\lambda_{1}} \pm1$ and
$\lambda_{1} \ge 3$. $\QED$


\section{Proof of the Alternating Sign Property}\label{alternatingsignprop}

Recall that the Alternating Sign Property (\rm{ASP}) is the
assertion that for any partition $\lambda \vdash n$,
\[ \sign(\eta_{\lambda}) = (-1)^{n-\lambda_{1}}= (-1)^{\textnormal{\# cells under the first row of }\lambda}. \]

\begin{prop}\label{SC-reduction-I}
Let $\lambda \vdash n$. Suppose \rm{ASP} holds for partitions of
smaller size. Then

\noindent $\bullet$ if $\lambda_{1} \equiv
\lambda_{2}~~(\textnormal{mod}~~2)$, then
$|\eta_{\lambda}|=|h|\eta_{\lambda-\w{c}}|-|\eta_{\lambda-\w{h}}||$;

\noindent $\bullet$ if $\lambda_{1} \not \equiv \lambda_{2}
~~(\textnormal{mod}~~2)$, then
$|\eta_{\lambda}|=h|\eta_{\lambda-\w{c}}|+|\eta_{\lambda-\w{h}}|$.

\noindent Moreover, \rm{ASP}~ holds for $\lambda$ if $\lambda_{1}
\not \equiv \lambda_{2} ~~(\textnormal{mod}~~2)$.  If $\lambda_1\equiv \lambda_2\ \ (\rm mod \ 2),$ then 
\rm{ASP} for $\lambda$ is equivalent to $|\eta_{\lambda-\w{h}}|<h|\eta_{\lambda-\w{c}}|$.

\end{prop}

\proof In view of the main recurrence
$\eta_{\lambda}=(-1)^{h}\left((-1)^{\lambda_{1}}h
\eta_{\lambda-\w{c}} + \eta_{\lambda-\w{h}}\right)$, the absolute
value of $\eta_{\lambda}$ depends on the values and the signs of
$(-1)^{\lambda_{1}}\eta_{\lambda-\w{c}}$ and $\eta_{\lambda-\w{h}}$.
Since \rm{ASP} holds for $\lambda-\w{c}$ and
$\lambda-\w{h}$, we have
\begin{eqnarray}
\sign((-1)^{\lambda_{1}}\eta_{\lambda-\w{c}}) & = &
(-1)^{\lambda_{1} +
(\# \textnormal{cells under the first row of } \lambda-\w{c})} \nonumber \\
& = & (-1)^{\lambda_{1} + \lambda_{2} - 1 + (\# \textnormal{cells
under the first row of } \lambda-\w{h})} \nonumber \\
& = & (-1)^{\lambda_{1}+\lambda_{2}-1} \sign(\eta_{\lambda-\w{h}}).
\label{sign}
\end{eqnarray}
Therefore, if $\lambda_{1} \not \equiv
\lambda_{2}~~(\textnormal{mod}~~2)$, then
$\sign((-1)^{\lambda_{1}}\eta_{\lambda-\w{c}})=\sign(\eta_{\lambda-\w{h}})$.
This implies that $\sign(\eta_{\lambda}) =
\sign((-1)^{h}(-1)^{\lambda_{1}}\eta_{\lambda-\w{c}}) =
(-1)^{\lambda_{1}+r-1+\lambda_{1}+(n-r)-(\lambda_{1}-1)}=
(-1)^{n-\lambda_{1}}$, i.e. \rm{ASP} holds for $\lambda$ and
$|\eta_{\lambda}|=h|\eta_{\lambda-\w{c}}|+|\eta_{\lambda-\w{h}}|$.
Otherwise,
$|\eta_{\lambda}|=|h|\eta_{\lambda-\w{c}}|-|\eta_{\lambda-\w{h}}||$.  Here \rm{ASP} is 
equivalent to $\sign(\eta_\lambda)=(-1)^{n-\lambda_1}$ which is equivalent 
to $|\eta_{\lambda-\w{h}}|<h|\eta_{\lambda-\w{c}}|$.
$\QED$

\begin{prop}\label{SC-reduction-II}
Let $\lambda \vdash n$, $n \ge 3$.  Assume $\lambda=(n)$ or 
$\lambda_{1} \ge \lambda_{2}+2$. Suppose \rm{ASP} holds for
partitions of smaller size. Then
\[ |\eta_{\lambda}| > |\eta_{(\lambda_{1}-2, \lambda_{2}, \ldots,
\lambda_{r})}|. \] Moreover, \rm{ASP} holds for $\lambda$.
\end{prop}

\proof Let $\lambda' = (\lambda'_{1}, \lambda_{2}, \ldots,
\lambda_{r})$ where $\lambda'_{1}=\lambda_{1}-2$, i.e. $\lambda'$ is
the partition of $n-2$ obtained from $\lambda$ by deleting the first
two cells of $\lambda_{1}$ from the right. Let $\w{h'}$ and $\w{c'}$
denote the hook and the first column of $\lambda'$ respectively.
Also, let $h'$ denote the size of $\w{h'}$. We shall prove by
induction on $n=|\lambda|$ that
$|\eta_{\lambda}|>|\eta_{\lambda'}|$.

When $n=3$, the only partition which satisfies the conditions of the
theorem is $\lambda = (3)$.  So $|\eta_{\lambda}|=2 >
|\eta_{\lambda'}|=0$. Indeed if $\lambda=(n)$ the statement follows as $D_n>D_{n-2}$.  This
means we can assume $r\ge 2$.
Let $n>3$.  As \rm{ASP} holds for
$\lambda-\w{c}$, $\lambda'-\w{c'}$, $\lambda-\w{h}$,
$\lambda'-\w{h'}$, we have
\begin{eqnarray}
\sign((-1)^{\lambda_{1}}\eta_{\lambda-\w{c}}) & = &
\sign((-1)^{\lambda_{1}-2}\eta_{\lambda-\w{c}}) \nonumber \\
& = & \sign((-1)^{\lambda'_{1}}\eta_{\lambda'-\w{c'}})
\nonumber \\
\sign(\eta_{\lambda-\w{h}}) & = & \sign(\eta_{\lambda'-\w{h'}}).
\nonumber
\end{eqnarray}
This means the signs in the main recurrence for $\eta_\lambda$ and $\eta_{\lambda'}$ are
the same.  In particular, $|\eta_{\lambda}| = |h|\eta_{\lambda-\w{c}}| +
|\eta_{\lambda-\w{h}}||$ if and only if $|\eta_{\lambda'}| =
|h'|\eta_{\lambda'-\w{c'}}| +  |\eta_{\lambda'-\w{h'}}||$ and the same when the signs
are negative.

We can use induction unless $\lambda=(3,1^{n-3})$.  By
Lemma~\ref{hook}, $|\eta_{(3,1^{n-3})}|=D_3+(n-3)D_2=n-1.$  However,
$|\eta_{(1^{n-2})}|=D_1+n-2-1=n-3$ and so the result holds.

We can now use induction, and so $|\eta_{\lambda-\w{c}}|>|\eta_{\lambda'-\w{c'}}|$. Since
$h=h'+2$ and $\lambda-\w{h} = \lambda'-\w{h'}$, we deduce that
\begin{eqnarray}
|\eta_{\lambda}| & = & |(h'+2)|\eta_{\lambda-\w{c}}| \pm |\eta_{\lambda'-\w{h'}}|| \nonumber \\
& > & |h'|\eta_{\lambda'-\w{c'}}| \pm
|\eta_{\lambda'-\w{h'}}|| \nonumber \\
& = & |\eta_{\lambda'}|. \nonumber
\end{eqnarray}
Here in the case with negative sign we have used $|\eta_{\lambda'-\w{h'}}|< h'|\eta_{\lambda'-\w{c'}}|$ 
since \rm{ASP} holds for $\lambda'$ using Proposition~\ref{SC-reduction-I}.

If $|\eta_{\lambda}|=h|\eta_{\lambda-\w{c}}|+|\eta_{\lambda-\w{h}}|$
then
$\sign(\eta_{\lambda})=\sign((-1)^{h}(-1)^{\lambda_{1}}\eta_{\lambda-\w{c}})=(-1)^{n-\lambda_{1}}$.
This follows as by assumption $\sign(\eta_{\lambda-\w
{c}})=(-1)^{n-r-(\lambda_1-1)}$ and so mod $2$, $\lambda_1+(r-1)
+\lambda_1+n-r-(\lambda_1-1)$ is $n-\lambda_1$.  This means \rm{ASP}
holds for $\eta_{\lambda}$.  Otherwise,
$|\eta_{\lambda}|=|h|\eta_{\lambda-\w{c}}|-|\eta_{\lambda-\w{h}}||$
and
$|\eta_{\lambda'}|=|h'|\eta_{\lambda'-\w{c'}}|-|\eta_{\lambda'-\w{h'}}||$.
But, since \rm{ASP} holds for $\lambda'$, this means that
$|\eta_{\lambda
'}|=h'|\eta_{\lambda'-\w{c'}}|-|\eta_{\lambda'-\w{h'}}| > 0$.
Therefore,
\begin{eqnarray}
h|\eta_{\lambda-\w{c}}| > h'|\eta_{\lambda'-\w{c'}}| \ge
|\eta_{\lambda'-\w{h'}}| = |\eta_{\lambda-\w{h}}|. \nonumber
\end{eqnarray}
This means $\sign(\eta_\lambda)$ is $(-1)^{h}(-1)^{\lambda_1}\sign(\eta_{\lambda-\w{c}})=(-1)^{n-\lambda_1}$
and so \rm{ASP} holds for $\eta_{\lambda}$. 
$\QED$

We want to prove \rm{ASP} for $\lambda$ by assuming that the
assertion holds for partitions of smaller size than $|\lambda|$. By
Proposition \ref{SC-reduction-I} and Proposition
\ref{SC-reduction-II}, it remains to consider the case
$\lambda_{1}=\lambda_{2}$.

To state our next results, we require some new terminology. For a
partition $\lambda=(\lambda_{1}, \ldots, \lambda_{r}) \vdash n$ and
$0 \le i \le \lambda_{1}$, let $\lambda-\w{c}_{i}$ denote the
partition obtained from $\lambda$ by deleting the first $i$ columns.
In particular, $\lambda = \lambda-\w{c}_{0}$, $\lambda-\w{c} =
\lambda-\w{c}_{1}$. Similarly, for $0 \le i \le r$, let $\lambda -
\w{\rho}_{i}$ denote the partition obtained from $\lambda$ by
deleting the first $i$ rows. When $i=1$, we also write $\lambda -
\w{\rho}$ instead of $\lambda-\w{\rho}_{1}$. Using these notations,
note that $\lambda - \w{h} = (\lambda - \w{c})-\w{\rho}$.

For the rest of this section, we let $h_{i}$ denote the size of the
hook of $\lambda - \w{c}_{i-1}$, where $1 \le i \le \lambda_{1}$. We
have the following upper bound for $|\eta_{\lambda}|$ in terms of
the $h_{i}$'s:

\begin{prop}\label{hooks}
Let $\lambda \vdash n$. Then \[ |\eta_{\lambda}| \le
\prod_{i=1}^{\lambda_{1}} (h_{i}+1). \]
\end{prop}

\proof As this is not needed in the sequel, we omit the proof which 
is a routine interation of Lemma~ $\QED$

We shall be interested in partitions with
$\lambda_{1}=\lambda_{2}$ and $\lambda_3<\lambda_1$ if $r\ge 3$ where as usual 
$r$ is the number of rows. For this we denote $\lambda_{1}$ by
$t$ and assume $t\geq 2$.  Note that the smallest partition satisfying these conditions is
$\lambda=(1^{2})$.  For the definition of $\delta$
below we assume $\lambda_1=\lambda_2$ with $\lambda_3<\lambda_1$ if $r\ge 3$.  We define the following functions:
\begin{eqnarray*}
H(\lambda) & = & \prod_{i=1}^{t-1}h_{i}-\prod_{i=1}^{t-2}h_{i} -
\prod_{i=2}^{t-1}(h_{i}-2)-\sum_{i=1}^{t-3} h_{1}h_{2} \cdots
h_{i}(h_{i+2}-2)(h_{i+3}-2)\cdots(h_{t-1}-2), \\
S(\lambda) & = & \prod_{i=1}^{t-1} (h_{i}-2), \nonumber
\\
\delta((1^{2})) & = & 1, \nonumber \\
\delta(\lambda) & = &
h_{1}\delta(\lambda-\w{c})-|\eta_{(\lambda-\w{c})-\w{\rho}_{2}}|
             \ \  \ {\rm where\  for\  this\  we\  assume} \ \lambda_1=\lambda_2\ {\rm with}\ \lambda_3<\lambda_1 \ {\rm 
if}\  r\ge 3.
\nonumber
\end{eqnarray*}
By convention we mean $H(\lambda)=h_{1}$ when $t=2$ and
$H(\lambda)=h_{1}h_{2}-h_{1}-(h_{2}-2)$ when $t=3$. Recursively, we
have
\begin{eqnarray}
\delta(\lambda) & = &
\prod_{i=1}^{t-1}h_{i}-\prod_{i=1}^{t-2}h_{i}-|\eta_{(\lambda-\w{c})-\w{\rho}_{2}}|-\sum_{i=1}^{t-3}
h_{1}h_{2} \cdots h_{i} |\eta_{\lambda-\w{c}_{i+1}-\w{\rho}_{2}}|.
\nonumber
\end{eqnarray}
The motivation for the above functions will become apparent in
Lemma~\ref{delta-lambda} and Proposition~\ref{lambda-lambda-row}.

\begin{lm}\label{S-lambda}
Let $\lambda \vdash n$, $\lambda_{1}=\lambda_{2} \ge 3$ and 
$\lambda_{3} <\lambda_{1}$ if $r\ge 3$. Then
$H(\lambda) > S(\lambda) > 0$.
\end{lm}

\proof We set $t= \lambda_{1} \ge 3$. Note that $h_{1}>h_{2}> \cdots
> h_{t} =2$. Clearly, $S(\lambda)>0$ as $h_{i} \ge 3$ for all $i \le
t-1$. We proceed by induction on $t$. For $t=3$, $h_{1}>h_{2} \ge
3$, and $H(\lambda)=h_{1}h_{2}-h_{1}-(h_{2}-2) > (h_{1}-2)(h_{2}-2)
= S(\lambda)>0$ since $h_{1}$ and $h_{2}$ are greater than $ 2$. Let $t>3$. Then
\begin{eqnarray}
H(\lambda) & = &
h_{1}\left(\prod_{i=2}^{t-1}h_{i}-\prod_{i=2}^{t-2}h_{i}-\prod_{i=3}^{t-1}(h_{i}-2)-\sum_{i=2}^{t-3}
h_{2} \cdots h_{i} (h_{i+2}-2) \cdots (h_{t-1}-2)\right) -
\prod_{i=2}^{t-1}(h_{i}-2). \nonumber
\end{eqnarray}
By the definition of $H(\lambda-\w{c})$ and the inductive
hypothesis,
\[ \left(\prod_{i=2}^{t-1}h_{i}-\prod_{i=2}^{t-2}h_{i}-\prod_{i=3}^{t-1}(h_{i}-2)-\sum_{i=2}^{t-3}
h_{2} \cdots h_{i} (h_{i+2}-2) \cdots (h_{t-1}-2)\right) =
H(\lambda-\w{c}) > S(\lambda-\w{c}). \] Therefore,
\begin{eqnarray}
H(\lambda) & > & h_{1}S(\lambda - \w{c}) -
\prod_{i=2}^{t-1}(h_{i}-2)
\nonumber \\
& = & h_{1}\prod_{i=2}^{t-1}(h_{i}-2) - \prod_{i=2}^{t-1}(h_{i}-2)
\nonumber \\
& = & (h_{1}-1)\prod_{i=2}^{t-1}(h_{i}-2) \nonumber \\
& > & \prod_{i=1}^{t-1}(h_{i}-2) = S(\lambda). \nonumber
\end{eqnarray} $\QED$

\begin{lm}\label{delta-lambda}
Let $\lambda \vdash n$, $\lambda_{1}=\lambda_{2}$ and $\lambda_{3} < \lambda_1$ if $r\ge 3$. Then $\delta(\lambda)
> 0$.
\end{lm}

\proof Since $\delta((1^{2}))=1>0$, we may assume that $\lambda \not
= (1^{2})$. If $\lambda_{1}=\lambda_{2}=2$ then $\delta(\lambda) =
h_{1}\delta((1^{2}))-1 >0$. So we may assume that
$\lambda_{1}=\lambda_{2} \ge 3$. By Lemma \ref{S-lambda}, it
suffices to show that $\delta(\lambda) \ge H(\lambda)$. Indeed, by
Proposition \ref{hooks}, we have, for all $i \in \{0, 1, \ldots,
t-3\}$,
\[ |\eta_{(\lambda-\w{c}_{i+1})-\w{\rho}_{2}}| \le \prod_{j=i+2}^{t-1}
(h_{j}-2).\] So
\begin{eqnarray}
\delta(\lambda) & \ge &
\prod_{i=1}^{t-1}h_{i}-\prod_{i=1}^{t-2}h_{i} -
\prod_{i=2}^{t-1}(h_{i}-2)-\sum_{i=1}^{t-3} h_{1}h_{2} \cdots
h_{i}(h_{i+2}-2)(h_{i+3}-2)\cdots(h_{t-1}-2) \nonumber \\
& = & H(\lambda). \nonumber
\end{eqnarray} $\QED$

\begin{prop}\label{lambda-lambda-row}
Let $\lambda \vdash n$, $\lambda_{1}=\lambda_{2}$ and $\lambda_{3}<\lambda_1$ if $r\ge 3$. Then
\[ |\eta_{\lambda}| \ge |\eta_{\lambda-\w{\rho}}| + \delta(\lambda) > |\eta_{\lambda-\w{\rho}}|.\]
\end{prop}

\proof The second inequality follows immediately from Lemma
\ref{delta-lambda}. We shall prove the first inequality by induction
on $\lambda_{1} \ge 1$. If $\lambda_{1}=1$, then $\lambda=(1^{2})$,
$|\eta_{\lambda}|=1$, $|\eta_{\lambda-\w{\rho}}|=0$ and
$\delta(\lambda)=1$, so the inequality holds. Let $\lambda_{1}>1$.
By induction, $|\eta_{\lambda-\w{c}}|
\geq |\eta_{(\lambda-\w{c})-\w{\rho}}|+\delta(\lambda-\w{c})$. Then using
Proposition~\ref{SC-reduction-I},
\begin{eqnarray}
|\eta_{\lambda}| & = &
h_{1}|\eta_{\lambda-\w{c}}|-|\eta_{(\lambda-\w{c})-\w{\rho}}| \nonumber \\
& \ge &
h_{1}|\eta_{\lambda-\w{c}}|-(|\eta_{\lambda-\w{c}}|-\delta(\lambda-\w{c}))~~(\textnormal{by
induction}) \nonumber \\
& = & (h_{1}-1)|\eta_{\lambda-\w{c}}|+\delta(\lambda-\w{c}).
\label{e7}
\end{eqnarray}
On the other hand,
\begin{eqnarray}
|\eta_{\lambda-\w{\rho}}| & \le &
(h_{1}-1)|\eta_{(\lambda-\w{c})-\w{\rho}}| +
|\eta_{(\lambda-\w{c})-\w{\rho}_{2}}| \nonumber \\
& \le &
(h_{1}-1)(|\eta_{\lambda-\w{c}}|-\delta(\lambda-\w{c}))+|\eta_{(\lambda-\w{c})-\w{\rho}_{2}}|
~~(\textnormal{by
induction})\nonumber \\
& = & (h_{1}-1)|\eta_{\lambda-\w{c}}| -
(h_{1}-1)\delta(\lambda-\w{c}) +
|\eta_{(\lambda-\w{c})-\w{\rho}_{2}}|. \label{e8}
\end{eqnarray}
As
$\delta(\lambda)=h_{1}\delta(\lambda-\w{c})-|\eta_{(\lambda-\w{c})-\w{\rho}_{2}}|$,
it follows from (\ref{e7}) and (\ref{e8}) that
\[ |\eta_{\lambda}| \ge |\eta_{\lambda-\w{\rho}}|+\delta(\lambda). \]
$\QED$

\begin{prop}\label{SC-reduction-III}
Let $\lambda \vdash n$, $\lambda_{1}=\lambda_{2}$. Suppose \rm{ASP}
holds for all partitions of size smaller than $n$. Then
\[ |\eta_{\lambda}| > |\eta_{\lambda-\w{\rho}}|.\]
\end{prop}

\proof If $\lambda_{3} <\lambda_1$ for $r\ge 3$ or $r=2$, 
then the assertion is true by Proposition \ref{lambda-lambda-row},
even without the assumptions that \rm{ASP} holds for smaller
partitions. So we may assume that $\lambda_{3} =\lambda_1$. We proceed by induction on
$\lambda_{1} \ge 1$. For $\lambda_{1}=1$, clearly $|\eta_{(1^{n})}|
= n-1 > n-2 = |\eta_{(1^{n-1})}|$. Let $\lambda_{1}>1$. By
induction, we assume that $|\eta_{\lambda-\w{c}}| >
|\eta_{(\lambda-\w{c})-\w{\rho}}|$. It follows from Proposition
\ref{SC-reduction-I} that
\begin{eqnarray}
|\eta_{\lambda}| & = &
h|\eta_{\lambda-\w{c}}|-|\eta_{(\lambda-\w{c})-\w{\rho}}|
\nonumber \\
& > & h|\eta_{\lambda-\w{c}}|-|\eta_{\lambda-\w{c}}| \nonumber \\
& = & (h-1)|\eta_{\lambda-\w{c}}| \nonumber \\
& > & (h-1)|\eta_{(\lambda-\w{c})-\w{\rho}}|. \label{e9}
\end{eqnarray}
On the other hand, since \rm{ASP} holds for
$\eta_{\lambda-\w{\rho}}$, $\eta_{(\lambda-\w{c})-\w{\rho}}$, 
and $\eta_{(\lambda-\w{c})-\w{\rho}_{2}}$, $\lambda_{2} =
\lambda_{3}$, and Proposition\ref{SC-reduction-I}
\begin{eqnarray}
|\eta_{\lambda-\w{\rho}}| & = &
h'|\eta_{(\lambda-\w{c})-\w{\rho}}|-|\eta_{(\lambda-\w{c})-\w{\rho}_{2}}| \nonumber \\
& < & h'|\eta_{(\lambda-\w{c})-\w{\rho}}| \nonumber \\
& = & (h-1)|\eta_{(\lambda-\w{c})-\w{\rho}}|, \label{e10}
\end{eqnarray}
where $h'=h-1$ is the hook of $\lambda-\w{\rho}$. It follows
immediately from (\ref{e9}) and (\ref{e10}) that
$|\eta_{\lambda}|>|\eta_{\lambda-\w{\rho}}|$. $\QED$

\begin{rem} \rm Notice it is not generally true that $|\eta_\lambda|>|\eta_{\lambda -\w{\rho}}|$ as if 
$\lambda=(2,1^{n-2})$, $|\eta_\lambda|=1$ and $|\eta_{\lambda -\w{\rho}}|=n-3$ by Lemma~\ref{hook}.  
\end{rem}

\noindent {\bf Proof of Theorem \ref{ASP}.}~~ We assume by induction
on the size of $\lambda$ that \rm{ASP} holds for partitions of size
smaller than $n$. Moreover, by Proposition \ref{SC-reduction-I} and
Proposition \ref{SC-reduction-II}, we may assume that
$\lambda_{1}=\lambda_{2}$. Recall that $h=\lambda_{1}+r-1$, where
$r$ is the length of the first column of $\lambda$. Then, by
Proposition \ref{SC-reduction-III},
$|\eta_{\lambda-\w{c}}|>|\eta_{(\lambda-\w{c})-\w{\rho}}|$, and so
$h|\eta_{\lambda-\w{c}}|>|\eta_{\lambda-\w{c}-\w{\rho}}|$. In view
of the main recurrence
\begin{eqnarray}
\eta_{\lambda} & = &
(-1)^{h}\left((-1)^{\lambda_{1}}h\eta_{\lambda-\w{c}}+\eta_{(\lambda-\w{c})-\w{\rho}}
\right), \nonumber
\end{eqnarray}
we deduce that
\begin{eqnarray}
\sign(\eta_{\lambda}) & = &
\sign\left((-1)^{h}(-1)^{\lambda_{1}}\eta_{\lambda-\w{c}} \right)
\nonumber \\
& = & (-1)^{h+\lambda_{1}} \cdot
(-1)^{(n-r)-(\lambda_{1}-1)}~~(\textnormal{by induction})
\nonumber \\
& = & (-1)^{n-\lambda_{1}}. \nonumber
\end{eqnarray}
$\QED$

\begin{cor}\label{consequenceASP}
Propositions~\ref{SC-reduction-I}, \ref{SC-reduction-II}, \ref{SC-reduction-III} hold
without the restriction that \rm{ASP} holds for all partitions of smaller size.  Furthermore,
in the case $\lambda_1\equiv \lambda_2$ (mod $2$), $|\eta_\lambda|=h|\eta_{\lambda-\w{c}}|-
|\eta_{\lambda-\w{h}}|$.

\end{cor}

\proof The first part follows from Theorem~\ref{ASP} which states that \rm{ASP} holds for all
partitions.  The second follows as $h|\eta_{\lambda-\w{c}}|>
|\eta_{\lambda-\w{h}}|$ as shown in the proof of Theorem~\ref{ASP}.
$\QED$


\section{A lower bound for $|\eta_{\lambda}|$}\label{lowerboundeta}

In this section we prove the lower bound for $\eta_\lambda$ when
$\lambda_1 \ge \lfloor\frac{n}{2}\rfloor$. This will be the proof of
Theorem~\ref{Main2}~(i). We begin by giving a lower bound for
$|\eta_{(a,a)}|$.

\begin{lm}\label{firstlowerbound}
Let $n=2a \ge 2$. Then $|\eta_{(a,a)}| \ge 2D_{a}+D_{a-1}$.
\end{lm}

\proof The assertion holds for $a=1,2,3$. We may assume $a \ge 4$
and proceed by induction on $a$. By the main recurrence and
induction,
\begin{eqnarray}
|\eta_{(a,a)}| & \ge & (a+1)|\eta_{(a-1,a-1)}|-|\eta_{(a-1)}|
\nonumber \\
& \ge & (a+1)(2D_{a-1}+D_{a-2})-D_{a-1} \nonumber \\
& = & 2(a+1)D_{a-1}+(a+1)D_{a-2}-D_{a-1} \nonumber \\
& = & 2aD_{a-1}+(a+1)D_{a-2}+D_{a-1} \nonumber \\
& \ge & 2(D_{a}-1) + (a+1)D_{a-2}+D_{a-1} ~~\textnormal{(by Lemma \ref{recurrence1})}\nonumber \\
& \ge & 2D_{a}+D_{a-1} \nonumber
\end{eqnarray}
since $(a+1)D_{a-2}-2 \ge 0$ for all $a \ge 4$. $\QED$

We are now ready to prove the lower bound of Theorem~\ref{Main2}~(i).

\begin{prop}\label{lowerbound}  Suppose $\lambda \vdash n$ with its first part equal to $\lambda_1 \ge \lfloor \frac{n}{2} \rfloor$.  Then
$|\eta_\lambda |\ge
D_{\lambda_1}+(n-\lambda_1)D_{\lambda_1-1}=|\eta_{(\lambda_1,1^{n-\lambda_1})}|.$  This is
the lower bound needed for Theorem~\ref{Main2}~(i).
\end{prop}
\proof We use induction on $n \ge 1$. Since all the cases for $n \le
11$ can be done by inspection, we let $n \ge 12$. Furthermore, by
Lemma~\ref{hook}, the assertion holds with equality when $\lambda$
is a hook regardless of $\lambda_{1}$. We may therefore assume that
$\lambda$ is not a hook.

We first show that the assertion holds for the following special
cases:

\noindent {\bf Case I.}~~$n=2b+1$ and $\lambda = (b,b,1)$.

\noindent By the main recurrence,
\begin{eqnarray}
|\eta_{\lambda}| & \ge & (b+2)|\eta_{(b-1,b-1)}|-|\eta_{(b-1)}| \nonumber \\
& \ge & (b+2)(2D_{b-1}+D_{b-2})-D_{b-1} ~~\textnormal{(by Lemma \ref{firstlowerbound})}\nonumber \\
& = & 2bD_{b-1}+3D_{b-1}+(b+2)D_{b-2} \nonumber \\
& \ge & 2(D_{b}-1)+3D_{b-1}+(b+2)D_{b-2}~~\textnormal{(by Lemma \ref{recurrence1})} \nonumber \\
& \ge & D_{b}+bD_{b-1}-3+3D_{b-1}+(b+2)D_{b-2}~~\textnormal{(by Lemma \ref{recurrence1})} \nonumber \\
& = & D_{b}+(b+1)D_{b-1}+2D_{b-1}+(b+2)D_{b-2}-3 \nonumber \\
& \ge &  D_{b}+(b+1)D_{b-1}, \nonumber
\end{eqnarray}
as for  $b\geq 4$, $2D_{b-1}+(b+2)D_{b-2}-3\geq 0$.

\noindent {\bf Case II.}~~$n=2b+1$ and $\lambda = (b,c,d)$ with $d\ge 2$.

\noindent Since $d \ge 2$, by Theorem \ref{Renteln_1},
$|\eta_{\lambda-\w{h}}|=|\eta_{(c-1,d-1)}| \le
\frac{D_{b-1}}{b-2}=D_{b-2}+D_{b-3}$. It follows from the main
recurrence that
\begin{eqnarray}
|\eta_{(b,c,d)}| & \ge &
(b+2)|\eta_{(b-1,c-1,d-1)}|-(D_{b-2}+D_{b-3}) \nonumber \\
& \ge & (b+2)\left(D_{b-1}+((2b+1)-3-(b-1))D_{b-2}\right)-(D_{b-2}+D_{b-3}) ~~\textnormal{(by induction)}\nonumber \\
& = &bD_{b-1}+2D_{b-1}+(b+2)(b-1)D_{b-2}-(D_{b-2}+D_{b-3}) \nonumber
\\
& \ge &
D_{b}-1+2D_{b-1}+(b+2)(D_{b-1}-1)-(D_{b-2}+D_{b-3})~~\textnormal{(by
Lemma \ref{recurrence1})} \nonumber
\\
& = & D_{b}+(b+1)D_{b-1}+3D_{b-1}-(b+3)-D_{b-2}-D_{b-3} \nonumber \\
& \ge & D_{b}+(b+1)D_{b-1},
\end{eqnarray}
as for $b\geq 5$, $3D_{b-1}-(b+3)-D_{b-2}-D_{b-3}\geq 0$.

From now on in this section, we may assume that either $\lambda_{1} > \frac{n}{2}$
or $\lambda_{1}=\frac{n}{2}$ and $\lambda$ has at least $3$ rows or
$\lambda_{1}=\frac{n-1}{2}$ and $\lambda$ has at least $4$ rows. The case
of $\lambda_1=\frac{n}{2}$ and just two rows is Lemma~\ref{firstlowerbound}.

In particular, $|\lambda-\w{h}| \le \lambda_{1}-2$ and so $|\eta_{\lambda-\w{h}}| \le D_{\lambda_{1}-2}$.

Let $r$ be the number of rows of $\lambda$. Suppose $\lambda -\w{c}$
has $s$ rows for some $r \ge s\ge 2$. Note that by the main
recurrence
\begin{eqnarray*}
 |\eta_\lambda| &\ge & h|\eta _{\lambda-\w{c}}|-|\eta_{\lambda-\w{h}}| \\
                                &\geq& (\lambda_1+r-1)|\eta_{\lambda-\w{c}}| - D_{\lambda_{1}-2}\\
                                &=& \lambda_1|\eta_{\lambda-\w{c}}| +(r-1)|\eta_{\lambda-\w{c}}|-D_{\lambda_{1}-2}.
\end{eqnarray*}
We proceed by induction on the size of $\lambda$ to conclude
$|\eta_{\lambda-\w{c}}|\ge
D_{\lambda_1-1}+(n-r-(\lambda_1-1))D_{\lambda_1-2}$.
\begin{eqnarray*}
|\eta_\lambda| &\ge &\lambda_1\big(D_{\lambda_1-1}+(n-r-(\lambda_1-1))D_{\lambda_1-2}\big) \\
               & &    +\ \ \   (r-1)\big(D_{\lambda_1-1}+(n-r-(\lambda_1-1))D_{\lambda_1-2}\big)-D_{\lambda_{1}-2} \\
                 &=&(D_{\lambda_1}\pm 1)+(n-\lambda_1-(r-1))(D_{\lambda_1-1}\pm1)+(n-\lambda_1-(r-1))D_{\lambda_1-2} \\
                  && \ \ \ +  (r-1)\big(D_{\lambda_1-1}+(n-r-(\lambda_1-1))D_{\lambda_1-2}\big)-D_{\lambda_{1}-2} \\
                &=&D_{\lambda_1}+(n-\lambda_1)D_{\lambda_1-1} \pm 1 \pm(n-\lambda_1-(r-1)) \\
           &&\ \ \ \ \ \ +r(n-\lambda_1-(r-1))D_{\lambda_1-2} -D_{\lambda_{1}-2}
\end{eqnarray*}
It remains to show that $\pm 1 \pm (n-h)+r(n-h)D_{\lambda_1-2}
-D_{\lambda_{1}-2}$ is positive. Notice
\begin{eqnarray*}
\pm 1 \pm (n-h)+r(n-h)D_{\lambda_1-2} -D_{\lambda_{1}-2} &\ge &
               -1 -(n-h)+r(n-h)D_{\lambda_1-2} - D_{\lambda_1-2}  \\
                &=& -1-n+h+\big(r(n-h)-1\big)D_{\lambda_1-2}.
\end{eqnarray*}
Note that $r\ge 2$, $\lambda_{1} \ge 4$ and $n-h\ge 1$. So
$-1-n+h+\big(r(n-h)-1\big)D_{\lambda_1-2}\ge
-(n-h)-1+D_{\lambda_1-2}\ge -(\lambda_1-1)-1+D_{\lambda_1-2}$ which is positive for $\lambda_1= 6$ and 
larger.

Hence,
the result follows. $\QED$


\section{An upper bound for $|\eta_{\lambda}|$}\label{upperboundeta}

In this section, we give the upper bound needed for
Theorem~\ref{Main2}~(i) for $|\eta_{\lambda}|$ except for some cases
when $\lambda_1=\frac{n}{2}$ for $n$ even or
$\lambda_1=\frac{n-1}{2}$ for $n$ odd.

\begin{prop}\label{etaupperbound-2}
Let $\lambda \vdash n$, with $\lambda_{1} \ge \frac{n-1}{2}$, $n \ge
2$, and suppose that $\lambda \not = (\frac{n}{2}, \frac{n}{2})$
when $n$ is even or $\lambda \not = (\frac{n-1}{2}, c, d)$ for any
$c \le \frac{n-1}{2}$, $d \ge 1$ when $n$ is odd. Then
\begin{eqnarray}
|\eta_{\lambda}| \le (n-\lambda_{1}+1)D_{\lambda_{1}} +
D_{\lambda_{1}-1}. \label{upperbound02}
\end{eqnarray}
\end{prop}

\proof It is readily checked that the theorem holds for $n=2,
\ldots, 13$. We shall let $n \ge 14$ so that $\lambda_{1} \ge 7$.
For $r=1$, $|\eta_{\lambda}|=D_{n} \le D_{n} + D_{n-1}$, and so the
theorem is true. Let $r \ge 2$.

If $\lambda$ is a hook then
$|\eta_{\lambda}|=D_{\lambda_{1}}+(n-\lambda_{1})D_{\lambda_{1}-1}$
(Lemma~\ref{hook}), which is clearly less than the right hand side
of (\ref{upperbound02}). So we may assume that $\lambda$ is not a
hook.

If $\lambda_{1} > \frac{n}{2}$, then $|\lambda-\w{h}| \le
\lambda_{1}-2$. If $\lambda_{1}=\frac{n}{2}$ (when $n$ is even), we
still have $|\lambda-\w{h}| \le \lambda_{1}-2$ since we are
excluding the case $\lambda =(\frac{n}{2},\frac{n}{2})$ so that $r
\ge 3$. If $\lambda_{1}=\frac{n-1}{2}$ (when $n$ is odd), we do have
$|\lambda-\w{h}| \le \lambda_{1}-2$ since we are excluding the shapes
$(\frac{n-1}{2}, c, d)$ with $3$ rows so that $r \ge 4$. Therefore,
in all cases, $|\eta_{\lambda-\w{h}}| \le D_{\lambda_{1}-2}$.

To validate the inductive hypothesis for $\eta_{\lambda-\w{c}}$,
$\lambda - \w{c} \vdash n-r$, we need to check that $\lambda_{1}-1
\ge \frac{n-r-1}{2}$, $n-r \ge 2$ and that $\lambda - \w{c} \not =
(\frac{n-r}{2}, \frac{n-r}{2})$ when $n-r$ is even or $\lambda -
\w{c} \not = (\frac{n-r-1}{2}, c', d')$ for any $c' \le
\frac{n-r-1}{2}$, $d' \ge 1$ when $n-r$ is odd. Indeed,
$\lambda_{1}- 1 \ge \frac{n-1}{2}-1 = \frac{n-3}{2} \ge
\frac{n-r-1}{2}$ as $r \ge 2$. Also, $n\geq \lambda_1+r-1$ and since
$\lambda_1\geq 7$, $n-r \ge 2$. Suppose
$\lambda-\w{c}=\big(\frac{n-r}{2},\frac{n-r}{2}\big)$. Then
$\frac{n-r}{2}=\lambda_1-1\geq\frac{n-3}{2}$ so that $r$ is $2$ or
$3$.  If $r=2$, the hypothesis gives $\lambda=(a,b)$ with $a>b$ as
$(\frac{n}{2},\frac{n}{2})$ is excluded so
$\lambda-\w{c}\neq(\frac{n-2}{2},\frac{n-2}{2})$.  This means $r=3$
and $\lambda = (\frac{n-1}{2}, \frac{n-1}{2}, 1)$ which is excluded
by the hypothesis. If $\lambda - \w{c} = (\frac{n-r-1}{2}, c', d')$
for some $c' \le \frac{n-r-1}{2}$ and $d' \ge 1$, then
$\frac{n-r-1}{2}=\lambda_{1}-1 \ge \frac{n-3}{2}$ which implies that
$r \le 2$, which is not possible because $\lambda-\w{c}$ has $3$
rows. Therefore, by induction,
\[ |\eta_{\lambda-\w{c}}| \le ((n-r)-(\lambda_{1}-1)+1)D_{\lambda_{1}-1} +
D_{\lambda_{1}-2}. \] Consequently,
\begin{eqnarray}
|\eta_{\lambda}| & \le & (\lambda_{1}+r-1)|\eta_{\lambda - \w{c}}| + |\eta_{\lambda-\w{h}}| \nonumber \\
& \le & (\lambda_{1}+r-1)|\eta_{\lambda - \w{c}}| +
D_{\lambda_{1}-2} \nonumber \\
& \le &
(\lambda_{1}+r-1)\left((n-r)-(\lambda_{1}-1)+1)D_{\lambda_{1}-1} +
D_{\lambda_{1}-2}\right)+D_{\lambda_{1}-2} \nonumber \\
& = & (\lambda_{1}+r-1)\left((n-\lambda_{1}+1)-(r-1)
\right)D_{\lambda_{1}-1} + (\lambda_{1}+r)D_{\lambda_{1}-2}
\nonumber \\
& = & \lambda_{1}(n-\lambda_{1}+1)D_{\lambda_{1}-1} +
(r-1)(n-2\lambda_{1}-r+2)D_{\lambda_{1}-1} +
(\lambda_{1}+r)D_{\lambda_{1}-2} \nonumber \\
& = & (n-\lambda_{1}+1)(D_{\lambda_{1}} \pm 1) +
(r-1)(n-2\lambda_{1}-r+2)D_{\lambda_{1}-1} +
(\lambda_{1}-1)D_{\lambda_{1}-2} + (r+1)D_{\lambda_{1}-2}\nonumber \\
& = &  (n-\lambda_{1}+1)D_{\lambda_{1}} \pm (n-\lambda_{1}+1) +
(r-1)(n-2\lambda_{1}-r+2)D_{\lambda_{1}-1} + D_{\lambda_{1}-1} \pm 1
+ (r+1)D_{\lambda_{1}-2}. \nonumber
\end{eqnarray}
It sufices to show that \begin{eqnarray}(n-\lambda_{1}+2) +
(r-1)(n-2\lambda_{1}-r+2)D_{\lambda_{1}-1} + (r+1)D_{\lambda_{1}-2}
\le 0.  \label{non-positive-1}\end{eqnarray}

The hypothesis on $\lambda_1$ gives $n - 2
\lambda_{1} \le 1$. Suppose $n - 2 \lambda_{1} = 1$ which because of
the partitions excluded implies $r \ge
4$. Then (\ref{non-positive-1}) becomes
\[ (\lambda_{1}+3)-(r-1)(r-3)D_{\lambda_{1}-1} +
(r+1)D_{\lambda_{1}-2} \le 0 \] By using the estimation
$(\lambda_{1}-1)D_{\lambda_{1}-2}-1 \le D_{\lambda_{1}-1}$, it
suffices to show that
\[ (\lambda_{1}+3)+(r+1)D_{\lambda_{1}-2} \le
(r-1)(r-3)\left((\lambda_{1}-1)D_{\lambda_{1}-2}-1\right), \] that
is \begin{eqnarray} \lambda_{1}+3+(r-1)(r-3)\le
\left((r-1)(r-3)(\lambda_{1}-1)-(r+1)  \right) D_{\lambda_{1}-2}.
\label{non-positive-3}
\end{eqnarray}
Since $\lambda_{1} \ge 7$, $D_{\lambda_{1}-2} \ge D_{5}=44$.
Therefore, it is enough to show that
\[  \lambda_{1}+3+(r-1)(r-3)     \leq 44\left((r-1)(r-3)(\lambda_{1}-1)-(r+1) \right). \]  For some of the computations
below we used the computer algebra package Maple.  Now we have
\[  45r^{2}-136r+182  \leq (44r^{2}-176r+131)\lambda_{1}. \]
As $\lambda$ is not a hook, we have $r \le \lambda_{1}$, and so we
are done if
\[  45r^{2}-136r+182  \leq (44r^{2}-176r+131)r , \]
\[ 0\leq 44r^{3}-221r^{2}+267r-182, \]
which indeed holds for all $r \ge 4$ again using Maple as for example $44r^3-221r^2$ is positive 
for $r\geq 6$ and smaller values can be checked.

Suppose first that $r=2$. As we are excluding the shape
$\lambda=(\frac{n}{2}, \frac{n}{2})$, we must have $\lambda=(a,b)$
with $a>b$ and so $n-2\lambda_{1} \le -1$.  From
(\ref{non-positive-1}) it is enough to show that \begin{eqnarray}
(n-\lambda_{1}+2)-D_{\lambda_{1}-1}+3D_{\lambda_{1}-2} & \le & 0\ \ \ \ {\rm or}
\nonumber \\
(n-\lambda_{1}+2)+3D_{\lambda_{1}-2} & \le &D_{\lambda_1-1}=
(\lambda_{1}-1)D_{\lambda_{1}-2} \pm 1, \nonumber
\end{eqnarray}
which is true if
\[ (n-\lambda_{1}+2) \le (\lambda_{1}-4)D_{\lambda_{1}-2}-1. \]
As $n-\lambda_{1}+2 \le \lambda_{1}+1$, it suffices to check that
\[ 1+ \frac{6}{\lambda_{1}-4} = \frac{\lambda_{1}+2}{\lambda_{1}-4} \le D_{\lambda_{1}-2}, \]
which clearly holds for all $\lambda_{1} \ge 7$.  This shows we may assume $r\geq 3$.

From now on, we may assume that $n-2\lambda_{1} \le 0$. Again, from
(\ref{non-positive-1}), it is enough to show that
\begin{eqnarray} (n-\lambda_{1}+2) +
(r-1)(-r+2)(\lambda_{1}-2)(D_{\lambda_{1}-2}+D_{\lambda_{1}-3}) +
(r+1)D_{\lambda_{1}-2} \le 0 \label{non-positive-2}\end{eqnarray}
where we used $D_{\lambda_1-1}=(\lambda_1-2)(D_{\lambda_1-2}+D_{\lambda_1-3})$

As $n-2\lambda_{1} \le 0$ and $r\geq 3$, rearranging
(\ref{non-positive-2}) we need to show that
\[
(n-\lambda_{1}+2)+\left((r+1)-(r-1)(r-2)(\lambda_{1}-2)
\right)D_{\lambda_{1}-2} \le
(r-1)(r-2)(\lambda_{1}-2)D_{\lambda_{1}-3}.\] Since $\lambda_{1}-2
\ge \frac{n}{2}-2$, the coefficient of $D_{\lambda_{1}-2}$ is less
than or equal to
\[ (r+1)-(r-1)(r-2)\big(\frac{n-4}{2}\big)  \]
which is nonpositive for $r \ge 3$, $n > 10$. It now suffices to show that
\[ (n-\lambda_{1}+2) \le (r-1)(r-2)(\lambda_{1}-2)D_{\lambda_{1}-3},\]
which is true for $\lambda_{1} \ge 7$, $r \ge 3$. $\QED$

\begin{prop}\label{sqeezebelow-2}
Suppose $\lambda \vdash  n$ with $\lambda_1 \ge \frac{n-1}{2}$, $n
\ge 3$, and $\lambda \not = (\frac{n}{2}, \frac{n}{2})$ when $n$ is
even or $\lambda \not = (\frac{n-1}{2}, c, d)$ for any $c \le
\frac{n-1}{2}$, $d \ge 1$ when $n$ is odd. Then
$|\eta_{\lambda}|<|\eta_{(\lambda_1+1,1^{n-\lambda_1-1})}|=
D_{\lambda_1+1}+(n-\lambda_1-1)D_{\lambda_1}$.  In particular the upper bound in
Theorem~\ref{Main2}~(i) holds except for the partitions excluded in
the statement.
\end{prop}
\proof In view of Proposition~\ref{etaupperbound-2}, we only need to
show that
$$(n-\lambda_1+1)D_{\lambda_1} +D_{\lambda_1-1} < D_{\lambda_1+1} +(n-\lambda_1-1)D_{\lambda_1}.$$
Subtracting $(n-\lambda_1-1)D_{\lambda_1}$ leaves
\begin{eqnarray*} 2D_{\lambda_1}+D_{\lambda_1-1}& < &D_{\lambda_1+1}   \\
                                                &=&(\lambda_1+1)D_{\lambda_1} \pm 1.
\end{eqnarray*}
Next, subtracting $2D_{\lambda_1}$ yields
\begin{eqnarray*}
D_{\lambda_1-1}          & <  & (\lambda_1-1)D_{\lambda_1} \pm 1
 \end{eqnarray*}
It is sufficient to show that
\begin{eqnarray*}
D_{\lambda_1-1}          & <  & (\lambda_1-1)D_{\lambda_1} - 1 \\
& = & (\lambda_{1}-1)^{2}(D_{\lambda_{1}-1}+D_{\lambda_{1}-2})-1
\end{eqnarray*}
This is true for $\lambda_1\ge 3$ (when $n \ge 7$) while the result
can be verified separately for small $n$. $\QED$


\section{Partitions with few parts}\label{fewpartslambda}

In this section we give the upper bounds needed for Theorem~\ref{Main2} (i) not covered by
Proposition~\ref{sqeezebelow-2}.

\subsection{The cases $(a,b)$, $(b,b)$ and $(b,b,1)$}

\begin{lm}\label{fewparts}  For $a\ge b>1$, the following formulae
hold:
\begin{eqnarray*}
 \eta_{(a,b)}  &=&(-1)^{a+1}D_{b-1}-(a+1)\eta_{(a-1,b-1)}  \\
 \eta_{(a,b-1,1)}  &=& (-1)^{a+2}D_{b-2} + (a+2)n_{(a-1,b-2)}
\end{eqnarray*}
\end{lm}

\proof This is just an application of the main
recurrence.  Notice
$(-1)^{a+1}(-1)^a=-1$ as one of $a$ and $a+1$ is even with the other
being odd.  Similarly $(-1)^{a+2}(-1)^a=1$ as they are either both
even or both odd. $\QED$.

We give an explicit formula for $\eta_{(a,b)}$ which is not specifically needed in the
remainder of the paper but gives some indication of how the values could be computed.

\begin{lm}\label{justtwoparts}  Suppose $(a,b) \vdash n$ with $b>0$.  Then
\begin{eqnarray*}
\eta_{(a,b)}  &=& (-1)^{a+1}\big(D_{b-1}+(a+1)D_{b-2}+(a+1)aD_{b-3} +\cdots +(a+1)a(a-1)\dots(a-b+3)D_0 \big)\\
                & & \ \ \ \  +\ (-1)^b(a+1)a(a-1)\dots (a-b+3)(a-b+2)D_{a-b}.
\end {eqnarray*}
\end{lm}

As this in not needed in the sequel, we omit the proof which is a straightforward 
iteration of Lemma~\ref{fewparts}

Next, we find the upper bounds needed for Theorem~\ref{Main2}~(i) for $|\eta_{(b,b)}|$ and
$|\eta_{(b,b,1)}|$ respectively.

\begin{lm}\label{betterbb}  For all $b$, $|\eta_{(b,b)}|\leq D_{b+1}+(b-1)D_b$.
If $b\ge 4$, then
$|\eta_{(b,b)}|\leq D_{b+1}+(b-3)D_b$ with equality only for $b=4$.  
\end{lm}

\proof The fist inequality is implied by the second for $b\geq 4$
and the small cases can be done by inspection.

For the second inequality we use induction on $b$ and note it is equality for $b=4$. Suppose $b>4$.
\begin{eqnarray*}
|\eta_{(b,b)}|  &= & h|\eta_{(b-1,b-1)}|-|\eta_{(b-1)}| \ \ \ \ \ \ ({\rm by \  Corollary~\ref{consequenceASP}}) \\
                 &=& (b+1)|\eta_{(b-1,b-1)}|   - D_{b-1}      \\
                 &\leq & (b+1)\big(D_b+(b-4)D_{b-1}\big) - D_{b-1}    \\
                  &=&  D_{b+1} \pm 1+(b-4)bD_{b-1}+(b-4)D_{b-1} -D_{b-1}     \\
                 &=& D_{b+1} \pm 1 +(b-4)D_b \pm(b-4) +(b-5)D_{b-1}
\end{eqnarray*}
Subtracting $D_{b+1}+(b-3)D_b$ we need
 $$\pm 1 -D_b \pm(b-4) +(b-5)D_{b-1}\leq 0.$$
This is
$$\pm 1\pm(b-4)+(b-5)D_{b-1}\leq D_b=bD_{b-1}\pm 1$$ which
is
$$\pm 1 \pm(b-4)\leq 5D_{b-1}\pm 1.$$
Taking the worst case of signs $+$ on the left and $-$ on the right this is
true if $b-2\leq 5D_{b-1}$ which is true for all $b\geq 5$.
Notice this is a strict inequality for $b\geq 3$.  We have already noted 
the lemma holds for $b=4$.  However, the induction step does not
apply with $b=4$.
$\QED$

We now consider the case $\eta_{(b,b,1)}$.

\begin{lm}\label{bb1}  Suppose $b>1$.  Then
$|\eta_{(b,b,1)}|  \leq  D_{b+1}+bD_b$ with equality only for $b=2$.
\end{lm}
\proof This can be checked by hand for $b\leq 4$ and so we assume
$b\ge 5$.
\begin{eqnarray*}
|\eta_{(b,b,1)}|  &  =  & (b+2)|\eta_{(b-1,b-1)}| -D_{(b-1)}   \ \ \ \ \ ({\rm by \  Corollary~\ref{consequenceASP}}) \\
               &\leq& (b+1+1)\big(D_b+(b-4)D_{b-1}\big)-D_{b-1}   \ \ \ \ \ ({\rm by~ Lemma~\ref{betterbb}})    \\
              &=&  D_{b+1} \pm 1+D_b +(b+2)(b-4)D_{b-1} - D_{b-1}  \\
                 &=&D_{b+1} \pm 1+D_b +(b-2)bD_{b-1}  -9D_{b-1}  \\
                 &=& D_{b+1} \pm1+D_b  +(b-2)D_{b}\pm (b-2) -9D_{b-1}  \\
                  &=&  D_{b+1} +(b-1)D_b  \pm 1  -9D_{b-1} \pm (b-2)
\end{eqnarray*}
Subtracting $D_{b+1}+bD_b$ gives $-D_b \pm 1  -9D_{b-1} \pm (b-2)$, which we must show is less
that $0$.  Taking the worst case with both signs $+$ we need $-9D_{b-1} +b-1<D_{b}.$  This is true if $b-1 <
D_b$ which is certainly true for $b\ge 5$. $\QED$

\subsection{The case $\lambda=(b,c,d)$ with $b=\frac{n-1}{2}$ and $d\geq 2$}

\begin{lm}\label{bcd}
Let $n=2b+1 \ge 3$. If $\lambda = (b,c,d) \vdash n$ with
$b=\frac{n-1}{2}$, $c<b$, $d \ge 2$. Then
\[ |\eta_{\lambda}| \le D_{b+1} + bD_{b}. \]
\end{lm}

\proof The smallest case is $n= 7$ where it holds.
So let $n \ge 9$ so that $b \ge 4$. Note that $|\lambda - \w{c}|=n-3
\ge 3$, $b-1 = \frac{n-3}{2}$. Moreover, $\lambda - \w{c} \not =
(\frac{n-3}{2},\frac{n-3}{2})$ since $c<b$ and also
$\lambda - \w{c} \not = (\frac{n-4}{2}, c', d')$ since
$n-3$ is even. By Proposition
\ref{etaupperbound-2},
\[ |\eta_{\lambda - \w{c}}| \le bD_{b-1}+D_{b-2}. \]
On the other hand, since $|\lambda-\w{h}|=n-(b+2)=b-1$, we have
\[ |\eta_{\lambda-\w{h}}| \le D_{b-1}. \]
Therefore,
\begin{eqnarray}
|\eta_{\lambda}| & \le & h|\eta_{\lambda-\w{c}}| +
|\eta_{\lambda-\w{h}}|
\nonumber \\
& \le & (b+2)(bD_{b-1}+D_{b-2})+D_{b-1} \nonumber \\
& = & ((b+2)b+1)D_{b-1} + (b+2)D_{b-2} \nonumber \\
& = & (b+1)^{2}D_{b-1}+(b-1)D_{b-2}+3D_{b-2} \nonumber \\
& = & (b+1)bD_{b-1}+(b+1)D_{b-1} + D_{b-1}\pm 1 + 3D_{b-2} \nonumber
\\
& \le & (b+1)(D_{b}+1) + bD_{b-1}+2D_{b-1}+3D_{b-2}+1 \nonumber \\
& \le & D_{b+1}+1+(b+1) + D_{b}+1+2D_{b-1}+3D_{b-2}+1 \nonumber \\
& = & D_{b+1} + D_{b} + 2D_{b-1}+3D_{b-2} + (b+4) \nonumber \\
& \le & D_{b+1}+bD_{b}, \nonumber
\end{eqnarray}
the last inequality holds since for all $b \ge 4$,
\[ (b-1)D_{b} = (b-1)^{2}D_{b-1}+(b-1)^{2}D_{b-2} \ge 9D_{b-1}+9D_{b-2} \ge 2D_{b-1}+3D_{b-2}+(b+4).
\] $\QED$

Note the exceptions to Proposition~\ref{sqeezebelow-2} with $\lambda_1=\frac{n-1}{2}$ other than $\eta_{(b,b,1)}$ have $d\geq 2$.  

\begin{cor}  Theorem~\ref{Main2}\label{firstupperbound} (i) has been proved.
\end{cor}
\proof  This follows from Proposition~\ref{sqeezebelow-2} and Lemmas \ref{betterbb}, \ref{bb1}, and \ref{bcd}.
$\QED$


\section{A second upper bound for $|\eta_{\lambda}|$}\label{secondupperboundeta}

In this section we provide another upper bound for $\eta_{\lambda}$
when $\lambda_1$ is small compared to $n$.  This is the upper bound needed for
Theorem~\ref{Main2}~(ii).  In particular we prove
the following theorem.
\begin{prop}\label{secondupperbound}  Let $b=\lfloor \frac{n}{2}\rfloor$ and suppose $\lambda \vdash n$.
If $\lambda_1<b$, then $|\eta_\lambda| \leq D_{b+1}+(n-b-1)D_{b}$ with strict inequality for $n\geq 6$.
\end{prop}

\begin{rem}\label{valuesnew}  \rm One implication of the proposition and
Corollay~\ref{firstupperbound} is that
$|\eta_\lambda | \leq D_{c+1} +(n-c-1)D_c$ for any $c$ greater than
or equal the maximum of $b=\lfloor \frac{n}{2}\rfloor$ and $\lambda_1$ by the use of
Lemma~\ref{hookcasevalues}.
\end{rem}

\proof   Assume $\lambda_1 < b$ and in particular $\lambda_1\leq b-1$.

We will use induction on $|\lambda|$.  As
$\lambda_1<\frac{n}{2}$ there must be at least three rows and we can
assume $\lambda=(\lambda_1,\lambda_2, \cdots ,\lambda_r)$ with $r\ge
3$.

If $\lambda$ were a hook the result would follow from
Lemma~\ref{hookcasevalues}.

We know from the main recurrence that
$$| \eta_\lambda| \leq h|\eta_{\lambda-\w{c}}| + |\eta_{\lambda-\w{h}}|.$$

Notice $\lambda_1-1 \leq b-2$ and $\lambda-\w{c}
=(\lambda_1-1,\lambda_2-1,\cdots ,\lambda_s-1)\vdash n-r$ with
$s\leq r$.  Also $\lambda_2-1\leq b-2$ as $\lambda_2 \leq \lambda_1$
and $\lambda-\w{h}=(\lambda_2-1, \cdots , \lambda_s-1)\vdash n-h$.  

We will apply the induction assumption to $\lambda-\w{c}$ and
$\lambda-\w{h}$.    Notice $\lfloor \frac{n-h}{2} \rfloor \leq b-2$
as $h\ge 4$. In particular $|\eta_{\lambda-\w{h}}|\leq D_{b-1}
+(n-h-b+1)D_{b-2}$ by Remark~\ref{valuesnew}.  Notice $n-h-b+1\leq b-1$ 
as $n-h\leq 2b-2$ as $\lfloor \frac{n-h}{2} \rfloor \leq b-2$.  
This means $|\eta_{\lambda -\w{h}} |\leq D_{b-1}+(b-1)D_{b-2}\leq D_{b-1}+D_{b-1}+1=2D_{b-1}+1$.

We now exclude the case in which $n$ is odd and $r=3$ which we will
do later.  With this assumption, $\lfloor \frac{n-r}{2} \rfloor \leq
b-2.$  Now by Remark~\ref{valuesnew} we see $|\eta_{\lambda-\w{c}}|\leq
D_{b-1} +(n-r-b+1)D_{b-2}$.

For these cases we have
\begin{eqnarray*}
|\eta_\lambda| &\leq & h|\eta_{\lambda -\w{c}}|+|\eta_{\lambda-\w{h}}| \\
               &\leq & h\big(D_{b-1}+(n-r-b+1)D_{b-2}\big) + |\eta_{\lambda -\w{h}}|.
\end{eqnarray*}
Notice $b=\lfloor \frac{n}{2}\rfloor$ and so $2b\geq n-1$.
Also notice that $h\leq n-1$ as $\lambda$ is not a hook and in
particular $hD_{b-1}\leq 2bD_{b-1}\leq 2D_b +2$. This means
$hD_{b-1} \leq 2D_{b}+2$.

Also $n-r-b+1\leq n-b-2\leq b-1$ as $b=\lfloor \frac{n}{2} \rfloor$.
In particular $(n-r-b+1)D_{b-2}\leq (b-1)D_{b-2}\leq D_{b- 1}+1$ and
so $h(n-r-b+1)D_{b-2}\leq 2b(D_{b-1}+1) \leq 2D_b+2+2b$.

In particular, 
\begin{eqnarray*}
|\eta_\lambda|  &\leq & 4D_b+2D_{b-1} +2b+5.
\end{eqnarray*}

But $4D_b+2D_{b-1} +2b+5 < D_{b+1}+(n-b-1)D_b$ using $D_{b+1}\ge
(b+1)D_{b}- 1$ as long as $b\ge 5$ and the smaller cases can be
done by hand.

We have one more case in which $n$ is odd, $r=3$ and $b'=\lfloor
\frac{n-3}{2}\rfloor=b-1$. Here $h=\lambda_1+2\leq b+1$.  By the
induction assumption and Remark~\ref{valuesnew},
$|\eta_{\lambda-\w{c}}|\leq D_b+(b-2)D_{b-1}$ as
$n-r-b=n-3-\frac{n-1}{2}=b-2$ and $\lambda_1-1<b'$. As shown above
$|\eta_{\lambda-\w{h}}|\leq 2D_{b-1}+1$.
\begin{eqnarray*}
|\eta_{\lambda}|   &\leq & (\lambda_1+2)|\eta_{\lambda -\w{c}}|+|\eta_{\lambda-\w{h}}|  \\
                  &\leq & (b +1)\big(D_b +(b-2)D_{b-1}\big)   +2D_{b-1}+1             \\
                    &= & D_{b+1} \pm 1 +  (b^2-b-2)D_{b-1}   +2D_{b-1}+1         \\
                     &\leq &  D_{b+1} + b^2D_{b-1} -bD_{b-1}-2D_{b-1} +2D_{b-1}+2               \\
                      &= & D_{b+1}+bD_b \pm b-D_b\pm 1+2.                             \\
                     &\leq &D_{b+1} + (b-1)D_b  + b+3.
\end{eqnarray*}

But this is less than $D_{b+1} +(n-b-1)D_b$ as long as $b\ge 4$ and
as usual the smaller cases follow by hand.  The needed condition is
$b+3 <D_{b}$ as here $n-b-1=b$. $\QED$

\begin{cor}  Theorem~\ref{Main2} has now been proven.
\end{cor}
\proof This follows from Corollary~\ref{firstupperbound} and Proposition~\ref{secondupperbound}.
$\QED$


\section{The case $\lambda_{1}=n-2$ }\label{casen-2}
In this section we prove the part of Theorem~\ref{Main3} when
$\lambda_1=n-2$. Before we begin, we need some preliminary
calculations.
\begin{lm}\label{P3}  Let $n\geq 4$.  $~~$
\begin{itemize}
\item[$(1)$] $\eta_{(n-2,2)} = \frac{(n-1)}{(n-3)}D_{n-2}$.
\item[$(2)$] $\eta_{(n-2,1^{2})} = nD_{n-3} + (-1)^{n}$.
\end{itemize}
\end{lm}
\proof The second is a direct application of Lemma~\ref{hook}.
The first follows from Lemma~\ref{fewparts} and Lemma~\ref{hook}.  $\QED$

\begin{lm}\label{P2}
Let $n \ge 6$. Then $|\eta_{(n-2,2)}|>|\eta_{(n-2,1^{2})}|$.
\end{lm}

\proof We must show that $\frac{n-1}{n-3}D_{n-2} =
\frac{n-1}{n-3}((n-3)(D_{n-3}+D_{n-4})) > nD_{n-3}+(-1)^{n}$, i.e.
$(n-1)(D_{n-3}+D_{n-4})> nD_{n-3}+(-1)^{n}$. Subtracting
$(n-1)D_{n-3}$ from both sides, we need to show that $(n-1)D_{n-4}>
D_{n-3} + (-1)^{n}$. Since $D_{n-3} = (n-3)D_{n-4} + (-1)^{n-3}$, it
suffices to show that $(n-1)D_{n-4}>(n-3)D_{n-4}$. But this is true
as $2D_{n-4}>0$ for all $n \ge 6$. $\QED$

\begin{cor}\label{casesn-2}  Theorem~\ref{Main3} has been proven when $\lambda_1=n-2$.  
\end{cor}

\proof This follows by Lemma~\ref{P2}.  $\QED$


\begin{section}{The case $\lambda_1=n-3$}\label{casen-3}

In this section we prove Theorem~\ref{Main3} in the case
$\lambda_1=n-3$.
\begin{lm}\label{PP1}  Let $n\geq 6$.
$~~$
\begin{itemize}
\item[$(1)$] $\eta_{(n-3,3)} = (-1)^{n-2}-\frac{(n-2)(n-3)}{(n-5)}D_{n-4}$.
\item[$(2)$] $\eta_{(n-3,2,1)} = -\frac{n-1}{n-4}D_{n-3}=-(n-1)(D_{n-4}+D_{n-5})$.
\item[$(3)$] $\eta_{(n-3,1^{3})}= (-1)^{n}(1+(-1)^{n-3}nD_{n-4})=-nD_{n-4}+(-1)^n$.
\end{itemize}
\end{lm}

\proof The third is Lemma~\ref{hook}.  The second follows from
Lemma~\ref{nearhook} and Lemma~\ref{recurrence1}.  For the first use the main recurrence and
substitute the value obtained in Lemma~\ref{P3}. $\QED$

\begin{lm}\label{P33}
Let $n \ge 6$. Then $|\eta_{(n-3,3)}|>|\eta_{(n-3,2,1)}|$.
\end{lm}

\proof We assume $n\geq 9$ and check the smaller values by hand.  We shall prove that
$\frac{n-1}{n-4}D_{n-3}<\frac{(n-2)(n-3)}{n-5}D_{n-4}-1$. As
$D_{n-3}-1 \le (n-3)D_{n-4}$, it is enough to show that
$\frac{n-1}{n-4}D_{n-3}<\frac{n-2}{n-5}D_{n-3}-\frac{n-2}{n-5}-1$,
which is $\frac{n-2}{n-5}+1 <
\left(\frac{n-2}{n-5}-\frac{n-1}{n-4}\right)D_{n-3} =
\frac{3}{(n-4)(n-5)}D_{n-3}$. As
$\frac{D_{n-3}}{n-4}=D_{n-4}+D_{n-5}$, it suffices to prove that
$\frac{2n-7}{3}<n-5<D_{n-4}+D_{n-5}$. This is true since
$D_{n-4}=(n-5)(D_{n-5}+D_{n-6})>n-5$. Recall $n\geq 9$ here.  $\QED$

\begin{lm}\label{P333}
Let $n \ge 6$. Then $|\eta_{(n-3,2,1)}|\geq |\eta_{(n-3,1^{3})}|$ with 
equality only for $n=6$.
\end{lm}

\proof Check this is equality if $n=6$ and so assume $n\geq 7$.  

We want to show that $nD_{n-4} + 1 < \frac{n-1}{n-4}D_{n-3}$.

\begin{eqnarray}
nD_{n-4}+1 & = & (n-3)D_{n-4}+3D_{n-4} + 1 \nonumber \\
           & \le & D_{n-3}+3D_{n-4}+2 \nonumber \\
           & < & \frac{n-1}{n-4}D_{n-3}. \nonumber
\end{eqnarray}
This will hold provided the following equivalent inequalities hold.  
\begin{eqnarray}
(n-4)D_{n-3}+2(n-4)+3(n-4)D_{n-4} & < & (n-1)D_{n-3} \nonumber \\
2(n-4)+3(n-4)D_{n-4} & < & 3D_{n-3} \nonumber \\
(n-4)(D_{n-4}+2/3) & < & D_{n-3}, \nonumber
\end{eqnarray}
which is true since $D_{n-3}=(n-4)(D_{n-4}+D_{n-5}) \ge
(n-4)(D_{n-4}+1) > (n-4)(D_{n-4}+2/3)$ for $n \ge 7$. $\QED$

\begin{cor}\label{casesn-3}  Theorem~\ref{Main3} is proven in the case $\lambda_1=n-3$.
\end{cor}
\proof This follow from Lemmas~\ref{P33} and \ref{P333}.
$\QED$
\end{section}


\begin{section}{The case $\lambda_1=n-4$}\label{casen-4}

In this section we prove Theorem~\ref{Main3} in the case
$\lambda_1=n-4$.  Notice $n\geq 8$ here.  

This will be proved by a series of lemmas.

\begin{lm}\label{valuesforanminus4}
For $a\ge 4$ the following values hold. \begin{eqnarray*}
\eta_{(a,4)}  &=&   2(-1)^{a+1}-(a+1)\eta_{(a-1,3)} \cr
\eta_{(a,3,1)} &=& (-1)^{a+2}
+(a+2)\big(D_{a-1}+2(D_{a-2}+D_{a-3})\big)     \cr 
\eta_{(a,2^2)}
&=&     (-1)^{a+1}(a+3) +(a+2)(a+1)D_{a-2} \cr \eta_{(a,2,1^2})  &=&
(a+3)\big((-1)^{a+3} + aD_{a-2}\big) \cr \eta_{(a,1^4)}  &=& (-1)^a
+ (a+4)D_{a-1}=D_a+4D_{a-1}.
\end{eqnarray*}
\end{lm}
\proof  These are straightforward applications of
the main recurrence, Lemma~\ref{recurrence1}, Lemma~\ref{hook}, Lemma~\ref{nearhook} and Lemma~\ref{P3} .
$\QED$

\begin{lm}\label{P4} If $a\ge 4$, then $\eta_{(a,4)}>\eta_{(a,3,1)}$.
\end{lm}

\proof By Lemma \ref{PP1} and Lemma \ref{valuesforanminus4},
\begin{eqnarray}
\eta_{(a,4)} & = &  2(-1)^{a+1}-(a+1)\eta_{(a-1,3)} \nonumber \\
& = &  2(-1)^{a+1}-(a+1)\left((-1)^{a} - \frac{a(a-1)}{(a-3)}D_{a-2}
\right) \nonumber \\
& \ge & \frac{(a+1)a(a-1)}{(a-3)}D_{a-2} - (a+3) \nonumber \\
& \ge & \frac{(a+1)a}{a-3}\left(D_{a-1}-1\right)-(a+3). \nonumber \\
& = & \frac{(a+1)a}{a-3}D_{a-1}-\left(\frac{(a+1)a}{a-3}+a+3\right)
\nonumber \\
& = & \frac{(a+1)a}{a-3}D_{a-1}-\left(2a+7+\frac{12}{a-3} \right).
\label{e-PP1-1}
\end{eqnarray}
On the other hand, by  Lemma \ref{valuesforanminus4},
\begin{eqnarray}
\eta_{(a,3,1)} &=& (-1)^{a+2}
+(a+2)\big(D_{a-1}+2(D_{a-2}+D_{a-3})\big) \nonumber  \\
& \le & 1+(a+2)\left(D_{a-1}+\frac{2}{a-2}D_{a-1} \right) \nonumber
\\
& = & 1+(a+2)\left(\frac{a}{a-2}D_{a-1} \right). \label{e-PP1-2}
\end{eqnarray}
Therefore, it is enough to show that the right hand side of
(\ref{e-PP1-1}) is more than the right hand side of (\ref{e-PP1-2}),
i.e.
\begin{eqnarray}
a\left(\frac{(a+1)(a-2)-(a+2)(a-3)}{(a-2)(a-3)}\right)D_{a-1} & >
& 2a+8+\frac{12}{a-3} \nonumber \\
\frac{4aD_{a-1}}{(a-2)(a-3)} & > & 2a+8+\frac{12}{a-3} \nonumber \\
\frac{4a}{a-3}(D_{a-2}+D_{a-3}) & > & 2a+8+\frac{12}{a-3} \nonumber \\
4a\left(D_{a-3}+D_{a-4}+\frac{D_{a-3}}{a-3}\right) & > &
2a+8+\frac{12}{a-3} \nonumber \\
D_{a-3}+D_{a-4}+\frac{D_{a-3}}{a-3} & > &
\frac{1}{2}+\frac{2}{a}+\frac{3}{(a-3)a} \label{e-PP1-3}
\end{eqnarray}
Note that the right hand side of (\ref{e-PP1-3}) is less than $2$
for $a \ge 4$. Therefore, the inequality holds whenever $a \ge 6$.
The lemma also holds for $a=4, 5$ by inspection.
 $\QED$

\begin{lm}\label{P44}  If $a\ge 4$, then $\eta_{(a,3,1)}>\eta_{(a,2^2)}$.
\end{lm}
\proof
\begin{eqnarray*}
\eta_{(a,3,1)}-\eta_{(a,2^2)}  &=&
(-1)^{a+2}-(-1)^{a+1}(a+3)+(a+2)\big(D_{a-1}+2(D_{a-2}+D_{a-3})-(a+1)D_{a-2}
\big) \cr
                             &=&(-1)^{a+2}(a+4)+(a+2)\big(D_{a-1}+2(D_{a-2}+D_{a-3})\big)-(a-1)D_{a-2}-2D_{a-2} \big)   \cr
              &=& (-1)^{a+2}(a+4) +(a+2)\big((a-1)D_{a-2} +(-1)^{a-1}+ 2(D_{a-2}+D_{a-3})           \cr
                     &&\ \ \ \ \ \ \ \ \ \ \ \ -(a-1)D_{a-2}-2D_{a-2} \big)   \cr
                     &=& (-1)^{a+2}(a+4) +(a+2)\big((-1)^{a-1}+2D_{a-3}\big)                         \cr
                 &=& (-1)^{a+2}(2) +(a+2)(2)D_{a-3}
\end{eqnarray*}
This is positive for $a=4$ and for $a\ge 5$ it is positive as $D_{a-3}\ge
1$. $\QED$

\begin{lm}\label{P444}  If $a\ge 4$, then $\eta_{(a,2,2)}>\eta_{(a,2,1^2)}$.
\end{lm}
\proof
\begin{eqnarray*}
\eta_{(a,2,2)}-\eta_{(a,2,1^2)} &=& (-1)^{a+1}(a+3)
+(a+2)(a+1)D_{a-2} -  (a+3)\big((-1)^{a+3} + aD_{a-2}\big)   \cr
                                &=& \big((a+2)(a+1)-(a+3)a\big)D_{a-2}             \cr
                                &=& (a^2+3a+2-a^2-3a)D_{a-2}                 \cr
                                &=& 2D_{a-2}.
\end{eqnarray*}
This is positive as $D_{a-2}\ge 1.$ $\QED$

\begin{lm}\label{P4444}  If $a\ge 4$, then $\eta_{(a,2,1^2)}>\eta_{(a,1^3)}$.
\end{lm}
\proof
\begin{eqnarray*}
\eta_{(a,2,1^2)}-\eta_{(a,1^4)} &=& (a+3)\big((-1)^{a+3} +
aD_{a-2}\big) -\big((-1)^a + (a+4)D_{a-1}\big)  \cr
                                &=& (a+4)(-1)^{a+3} +(a+3)aD_{a-2}-(a+4)\big((a-1)D_{a-2}+(-1)^{a-1}\big) \cr
                                &=& (a+4)\big((-1)^{a+3}-(-1)^{a-1}\big) + \big((a+3)a-(a+4)(a-1)\big)D_{a-2} \cr
                                &=& \big(a^2+3a-(a^2+3a-4)\big)D_{a-2}  \cr
                                &=& 4D_{a-2}.
\end{eqnarray*}
This is positive as $D_{a-2}\ge 1$.$\QED$

\begin{cor}\label{casesn-4}  Theorem~\ref{Main3} has been proven for the case $\lambda_1=n-4$.
\end{cor}
\proof This is Lemmas~\ref{P4}, \ref{P44}, \ref{P444}, and \ref{P4444}.
$\QED$

This completes the proof of Theorem~\ref{Main3} using Theorem~\ref{Renteln_1}, 
Corollary~\ref{casesn-2}, Corollary~\ref{casesn-3} and Corollary~\ref{casesn-4}.

\end{section}


\begin{section}{Some Values of $\eta_\lambda$}\label{values}

In this section we tabulate values of $\eta_\lambda$ for some small
values of $n$.

\medskip


{\small
\begin{center}
\begin{minipage}[t]{5cm}
\begin{tabular}{|rr|}
\multicolumn{2}{c} {$n=2$}  \\
\hline
$\lambda$ &  $\eta_{\lambda}$  \\
\hline
$2$& $1$  \\
$1^{2}$& $-1$  \\
\hline
\end{tabular}
\end{minipage}
\hspace{2cm}
\begin{minipage}[t]{5cm}
\begin{tabular}{|rr|}
\multicolumn{2}{c} {$n=3$} \\
\hline
$\lambda$ &  $\eta_{\lambda}$  \\
\hline
$3$& $2$ \\
$2,1$& $-1$ \\
$1^{3}$& $2$ \\
\hline
\end{tabular}
\end{minipage}
\end{center}


\begin{center}
\begin{minipage}[t]{5cm}
\begin{tabular}{|rr|r|rr|}
\multicolumn{5}{c} {$n=4$} \\
\hline
$\lambda$ &  $\eta_{\lambda}$ & \ \ \ & $\lambda$ &  $\eta_{\lambda}$  \\
\hline
$4$& $9$ & & $2,1^{2}$ & $1$ \\
$3,1$& $-3$ & &  $1^{4}$ & $-3$ \\
$2,2$& $3$ & &  &  \\
\hline
\end{tabular}
\end{minipage}
\hspace{2cm}
\begin{minipage}[t]{5cm}
\begin{tabular}{|rr|r|rr|}
\multicolumn{5}{c} {$n=5$} \\
\hline
$\lambda$ &  $\eta_{\lambda}$ & \ \ \ & $\lambda$ &  $\eta_{\lambda}$  \\
\hline
$5$& $44$ & & $2^{2},1$ & $-4$ \\
$4,1$& $-11$ & & $2,1^{3}$ & $-1$ \\
$3,2$& $4$ & & $1^{5}$ & $4$ \\
$3,1^{2}$& $4$ & &  &  \\
\hline
\end{tabular}
\end{minipage}
\end{center}


\begin{center}
\begin{minipage}[t]{5cm}
\begin{tabular}{|rr|r|rr|}
\multicolumn{5}{c} {$n=6$} \\
\hline
$\lambda$ &  $\eta_{\lambda}$ & \ \ \ & $\lambda$ &  $\eta_{\lambda}$  \\
\hline
$ 6 $  & $265$ & & $ 3, 1^{3}$ & $-5$\\
$ 5, 1 $  & $-53$ & & $ 2^{3} $ & $7$\\
$ 4, 2 $ & $15$ & & $ 2^{2}, 1^{2} $ & $5$\\
$4, 1^{2} $ & $13$ & & $ 2, 1^{4} $ & $ 1$\\
$3^{2}$ & $-11$ & & $ 1^{6} $ & $-5$ \\
$3, 2, 1$ & $-5$ & & &\\
\hline
\end{tabular}
\end{minipage}
\hspace{2cm}
\begin{minipage}[t]{5cm}
\begin{tabular}{|rr|r|rr|}
\multicolumn{5}{c} {$n=7$} \\
\hline
$\lambda$ &  $\eta_{\lambda}$ & \ \ \ & $\lambda$ &  $\eta_{\lambda}$  \\
\hline
$ 7 $ & $1854$  & & $3, 2^{2} $ & $ 6$\\
$6, 1 $ & $-309$ & & $3, 2, 1^{2}$ & $6$\\
$5, 2 $ & $66$ & & $3, 1^{4} $ & $ 6$\\
$5, 1, 1$ & $62$ & &$2^{3}, 1 $ & $-9$ \\
$4, 3 $ & $-21$ & &$2^{2}, 1^{3} $ & $-6$ \\
$4, 2, 1$ & $-18$ & &$2, 1^{5}$ & $ -1$ \\
$4, 1^{3} $ & $ -15$ & &$1^{7}$ & $ 6$ \\
$3^{2}, 1 $ & $ 14$ & & & \\
\hline
\end{tabular}
\end{minipage}
\end{center}


\begin{center}
\begin{minipage}[t]{5cm}
\begin{tabular}{|rr|r|rr|}
\multicolumn{5}{c} {$n=8$} \\
\hline
$\lambda$ &  $\eta_{\lambda}$ & \ \ \ & $\lambda$ &  $\eta_{\lambda}$  \\
\hline
$ 8 $ & $ 14833$ & & $4, 1^{4}$ & $ 17$\\
$7, 1$ &  $-2119$ &&$3^{2}, 2$ & $-19$ \\
$6, 2$ & $371$ & & $3^{2}, 1^{2}$ & $ -17$\\
$6, 1^{2}$ & $353$ & & $3, 2^{2}, 1$ & $-7$\\
$5, 3$ & $-89$ & & $3, 2, 1^{3}$ & $-7$\\
$5, 2, 1$ & $-77$ & & $3, 1^{5}$ & $-7$\\
$5, 1^{3}$ & $ -71$ & & $2^{4}$ & $13$\\
$4^{2}$ & $53$ & &$2^{3}, 1^{2}$ & $11$ \\
$4, 3, 1$ & $25$ & & $2^{2}, 1^{4}$ & $7$\\
$4, 2^{2}$ & $ 23$ & & $2, 1^{6}$ & $1$ \\
$4, 2, 1^{2}$ & $ 21$ & &$1^{8}$ & $-7$ \\
\hline
\end{tabular}
\end{minipage}
\hspace{2cm}
\begin{minipage}[t]{5cm}
\begin{tabular}{|rr|r|rr|}
\multicolumn{5}{c} {$n=9$} \\
\hline
$\lambda$ &  $\eta_{\lambda}$ & \ \ \ & $\lambda$ &  $\eta_{\lambda}$  \\
\hline
$9$ & $ 133496$ &&$4, 2^{2}, 1$ & $  -27$ \\
$8, 1$ & $-16687$ && $4, 2, 1^{3}$ & $  -24$\\
$7, 2$ & $2472$ &&$4, 1^{5}$ & $  -19$ \\
$7, 1^{2}$ & $2384$ && $3^{3}$ & $  32$\\
$6, 3$ & $ -463$ &&$3^{2}, 2, 1$ & $ 23$ \\
$6, 2, 1$ & $-424$ && $3^{2}, 1^{3}$ & $  20$\\
$6, 1^{3}$ & $-397$ &&$3, 2^{3}$ & $    8$ \\
$5, 4$ & $128$ && $3, 2^{2}, 1^{2}$ & $  8$\\
$5, 3, 1$ & $104$ && $3, 2, 1^{4}$ & $  8$\\
$5, 2^{2}$ & $92$ &&$3, 1^{6}$ & $  8$ \\
$5, 2, 1^{2}$ & $ 88$ &&$2^{4}, 1$ & $  -16$ \\
$5, 1^{4}$ & $  80$ && $2^{3}, 1^{3}$ & $  -13$\\
$4^{2}, 1$ & $  -64$ &&$2^{2}, 1^{5}$ & $  -8$ \\
$4, 3, 2$ & $  -31$ && $2, 1^{7}$ & $ -1$\\
$4, 3, 1^{2}$ & $  -29$ && $1^{9}$ & $  8$\\
\hline
\end{tabular}
\end{minipage}
\end{center}

\begin{center}
\begin{tabular}{|rr|r|rr|r|rr|r|rr|}
\multicolumn{11}{c} {$n=10$} \\
\hline
$\lambda$ &  $\eta_{\lambda}$ & \ \ \ & $\lambda$ &  $\eta_{\lambda}$ & \ \ \  &  $\lambda$ & $\eta_{\lambda}$&  \ \ \ &$\lambda$ & $\eta_{\lambda}$ \\
\hline
$ 10  $ & $  1334961 $ & & $  6,1^4  $ & $ 441  $ & & $ 4,3,2,1$ & $36$ & &$ 3,2^2,1^3$ & $-9$ \\
$ 9,1   $ & $ -148329  $ & & $5,5 $ & $-309 $  & & $4,3,1^3 $ & $33$& &$3,2,1^5$  &$-9$ \\
$  8,2  $ & $  19071 $ & &  $5,4,1  $ & $ -149$  & & $4,2^3 $ & $33$& &$3,1^7$  &$-9$ \\
$  8,1^2  $ & $  18541 $ & & $5,3,2 $ & $-125 $  & & $4,2^2,1^2 $ & $31$& &$2^5$  &$21$ \\
$  7,3  $ & $   -2967$ & & $5.3.1^2 $ & $-119 $  & & $4,2,1^4 $ & $27$& &$2^4,1^2$  &$19$ \\
$ 7,2,1   $ & $  -2781$ & & $5,2^2,1 $ & $-105 $  & & $4,1^6 $ & $21$& &$2^3,1^4$  &$15$ \\
$ 7,1^3   $ & $ -2649  $ & & $5,2,1^3 $ & $-99 $  & & $3^3,1 $ & $-39$& &$2^2,1^6$  &$9$ \\
$  6,4  $ & $  621 $ & & $5,1^5 $ & $-89 $  & & $3^2,2^2 $ & $-29$& & $2,1^8$ &$1$  \\
$ 6,3,1  $ & $529   $ & & $4^2,2 $ & $81 $  & & $3^2,2,1^2 $ & $-27$& & $1^{10}$ &$-9$ \\
$ 6,2^2   $ & $495   $ & & $4^2,1^2 $ & $75 $  & & $3^2,1^4 $ & $-23$& &  & \\
$ 6,2,1^2   $ & $  477 $ & & $4,3^2 $ & $39 $  & & $3,2^3,1 $ & $-9$& &  & \\
\hline
\end{tabular}
\end{center}
}

{\small

\begin{center}
\begin{tabular}{|rr|r|rr|r|rr|r|rr|}
\multicolumn{11}{c} {$n=11$,\  $\lambda_{1} \ge 5$} \\
\hline
$\lambda$ &  $\eta_{\lambda}$ & \ \ \ & $\lambda$ &  $\eta_{\lambda}$ & \ \ \  &  $\lambda$ & $\eta_{\lambda}$&  \ \ \ &$\lambda$ & $\eta_{\lambda}$ \\
\hline
$  11 $ & $ 14684570  $ & & $  7, 3, 1  $ & $ 3338   $  & & $ 6, 2^{2}, 1$ & $-557 $ & & $ 5, 3, 1^{3}$  & $134$ \\
$  10, 1 $ & $ -1468457  $ & & $ 7, 2^{2}  $ & $  3178   $  & & $6, 2, 1^{3} $ & $-530 $ & & $5, 2^{3} $  & $122$ \\
$  9, 2 $ & $ 166870  $ & & $ 7, 2, 1^{2}  $ & $ 3090   $  & & $6, 1^{5} $ & $-485 $ & & $5, 2^{2}, 1^{2} $  & $118$ \\
$ 9, 1^{2}  $ & $ 163162  $ & & $ 7, 1^{4}   $ & $  2914  $  & & $5^{2}, 1 $ & $362 $ & & $5, 2, 1^{4} $  & $ 110$ \\
$ 8, 3  $ & $  -22249  $ & & $ 6, 5  $ & $ -905   $  & & $5, 4, 2 $ & $ 178 $ & & $5, 1^{6} $  & $98$ \\
$ 8, 2, 1  $ & $  -21190 $ & & $ 6, 4, 1  $ & $ -710   $  & & $ 5, 4, 1^{2}$ & $170 $ & & $ $  & $$ \\
$ 8, 1^{3}  $ & $ -20395  $ & & $ 6, 3, 2  $ & $ -617   $  & & $ 5, 3^{2} $ & $158 $ & & $ $  & $$ \\
$ 7, 4  $ & $  3706  $ & & $ 6, 3, 1^{2}  $ & $ -595   $  & & $5, 3, 2, 1 $ & $143 $ & & $ $  & $$  \\

\hline
\end{tabular}
\end{center}

\begin{center}
\begin{tabular}{|rr|r|rr|r|rr|r|rr|}
\multicolumn{11}{c} {$n=12$,\  $\lambda_{1} \ge 6$} \\
\hline
$\lambda$ &  $\eta_{\lambda}$ & \ \ \ & $\lambda$ &  $\eta_{\lambda}$ & \ \ \  &  $\lambda$ & $\eta_{\lambda}$&  \ \ \ &$\lambda$ & $\eta_{\lambda}$ \\
\hline
$ 12  $ & $ 176214841  $ & & $ 8, 3, 1   $ & $ 24721   $  & & $ 7, 2^{2}, 1 $ & $ -3531$ & & $6, 3, 2, 1 $  & $694$ \\
$   11, 1 $ & $  -16019531  $ & & $ 8, 2^{2}   $ & $  23839   $  & & $7, 2, 1^{3}  $ & $ -3399$ & & $6, 3, 1^{3}  $  & $661$ \\
$ 10, 2  $ & $  1631619  $ & & $8, 2, 1^{2}  $ & $23309  $  & & $7, 1^{5} $ & $  -3179$ & & $ 6, 2^{3}$  & $637$ \\
$ 10, 1^{2}  $ & $ 1601953  $ & & $ 8, 1^{4}   $ & $ 22249   $  & & $ 6^{2} $ & $ 2119$ & & $6, 2^{2}, 1^{2} $  & $619 $ \\
$ 9, 3  $ & $ -190709   $ & & $7, 5 $ & $ -4959   $  & & $ 6, 5, 1$ & $ 1033 $ & & $6, 2, 1^{4}  $  & $583$ \\
$ 9, 2, 1  $ & $ -183557  $ & & $ 7, 4, 1  $ & $  -4169   $  & & $6, 4, 2 $ & $ 829 $ & & $6, 1^{6} $  & $ 529$ \\
$ 9, 1^{3}  $ & $ -177995  $ & & $7, 3, 2  $ & $  -3815  $  & & $6, 4, 1^{2}  $ & $ 799$ & & $ $  & $$ \\
$ 8, 4  $ & $  26701   $ & & $ 7, 3, 1^{2}  $ & $ -3709   $  & & $6, 3^{2} $ & $ 739$ & & $ $  & $$  \\

\hline
\end{tabular}
\end{center}

\begin{center}
\begin{tabular}{|rr|r|rr|r|rr|r|rr|}
\multicolumn{11}{c} {$n=13$, \ $\lambda_{1} \ge 6$} \\
\hline
$\lambda$ &  $\eta_{\lambda}$ & \ \ \ & $\lambda$ &  $\eta_{\lambda}$ & \ \ \  &  $\lambda$ & $\eta_{\lambda}$&  \ \ \ &$\lambda$ & $\eta_{\lambda}$ \\
\hline
$ 13  $ & $  2290792932  $ & & $  9, 1^{4}   $ & $ 192828   $  & & $7, 4, 1^{2} $ & $4632 $ & & $6, 4, 3 $  & $-996$ \\
$ 12, 1  $ & $ -190899411  $ & & $ 8, 5  $ & $ -33363    $  & &$7,3^{2} $ & $4452$ & & $6, 4, 2, 1 $  & $ -933$ \\
$ 11, 2  $ & $ 17621484  $ & & $ 8, 4, 1 $ & $ -29668  $  & & $7, 3, 2, 1 $ & $4239 $ & & $ 6, 4, 1^{3} $  & $-888$ \\
$ 11, 1^{2}  $ & $  17354492  $ & & $ 8, 3, 2   $ & $  -27811   $  & & $7, 3, 1^{3}  $ & $ 4080$ & & $6, 3^{2}, 1 $  & $ -831$ \\
$ 10, 3  $ & $ -1835571   $ & & $8, 3, 1^{2} $ & $   -27193  $  & & $7, 2^{3} $ & $ 3972  $ & & $6, 3, 2^{2} $  & $ -793$ \\
$ 10, 2, 1  $ & $ -1779948  $ & & $ 8, 2^{2}, 1  $ & $  -26223   $  & & $ 7, 2^{2}, 1^{2} $ & $3884 $ & & $6, 3, 2, 1^{2} $  & $-771$ \\
$ 10, 1^{3}  $ & $ -1735449  $ & & $ 8, 2, 1^{3} $ & $  -25428   $  & & $ 7, 2, 1^{4} $ & $3708 $ & & $6, 3, 1^{4} $  & $-727$ \\
$ 9, 4  $ & $  222492   $ & & $ 8, 1^{5}  $ & $  -24103   $  & & $7, 1^{6} $ & $3444 $ & & $ 6, 2^{3}, 1$  & $-708$  \\
$ 9, 3, 1  $ & $ 209780  $ & & $7, 6   $ & $ 7284   $  & & $6^{2}, 1  $ & $-2428 $ & & $6, 2^{2}, 1^{3} $  & $-681$ \\
$ 9, 2^{2}  $ & $  203952  $ & & $7, 5, 1  $ & $ 5580   $  & & $6, 5, 2  $ & $ -1203 $ & & $ 6, 2, 1^{5}$  & $ -636$ \\
$  9, 2, 1^{2} $ & $ 200244   $ & & $7, 4, 2   $ & $  4764  $  & & $6, 5, 1^{2} $ & $-1161 $ & & $6, 1^{7} $  & $-573$  \\

\hline
\end{tabular}
\end{center}
}

{\small

\begin{minipage}[t]{5cm}
\begin{tabular}{|rr|r|rr|r|rr|r|rr|}
\multicolumn{11}{c} {$n=15$} \\
\hline
$\lambda$ &  $\eta_{\lambda}$ & \ \ \ & $\lambda$ &  $\eta_{\lambda}$ & \ \ \  &  $\lambda$ & $\eta_{\lambda}$&  \ \ \ &$\lambda$ & $\eta_{\lambda}$ \\
\hline
$ 15  $ & $  481066515734 $ & & $  7^{2}, 1   $ & $ 18806  $ & & $ 6, 2, 1^{7} $ & $-742$ & &$4, 3^{2}, 2^{2}, 1$ & $-77$ \\
$ 14, 1 $ & $ -34361893981 $ & & $  7, 6, 2  $ & $ 9350 $  & & $ 6, 1^{9} $ & $-661$ & & $4, 3^{2}, 2, 1^{3}$ & $-74$ \\
$ 13, 2  $ & $  2672591754 $ & & $  7, 6, 1^{2} $ & $  9094 $  & & $ 5^{3} $ & $1214$ & & $4, 3^{2}, 1^{5}$ & $ -69$   \\
$ 13, 1^{2} $ & $ 2643222614 $ & & $ 7, 5, 3   $ & $ 7446 $ & & $5^{2}, 4, 1 $ & $859$ & & $4, 3, 2^{4}$ & $-73$ \\
$ 12, 3  $ & $ -229079293 $ & & $ 7, 5, 2, 1  $ & $  7089 $ & & $5^{2}, 3, 2 $ & $742$ & & $4, 3, 2^{3}, 1^{2}$ & $-71$  \\
$ 12, 2, 1  $ & $ -224273434 $ & & $ 7, 5, 1^{3}  $ & $   6822 $ & &$5^{2}, 3, 1^{2} $ & $714$ & & $ 4, 3, 2^{3}, 1^{4}$ & $-67$  \\
$ 12, 1^{3}  $ & $  -220268551 $ & & $  7, 4^{2} $ & $  6662 $ & & $5^{2}, 2^{2}, 1 $ & $662$& & $4, 3, 2, 1^{6}$ & $-61$  \\
$ 11, 4  $ & $  22026854  $ & &$ 7, 4, 3, 1  $ & $  6174 $ & &$5^{2},2,1^{3}  $ & $629$ & & $ 4, 3, 1^{8}$ & $-53$   \\
$ 11, 3, 1 $ & $ 21211046  $ & & $ 7, 4, 2^{2} $ & $  5954 $ & &$5^{2}, 1^{5}$ & $574$ & & $4, 2^{5}, 1 $ & $-66$  \\
$ 11, 2^{2} $ & $  20825390 $ & & $ 7, 4, 2, 1^{2} $ & $ 5822 $ & & $5, 4^{2}, 2$ & $374$& & $4, 2^{4}, 1^{3}$ & $-63$  \\
$ 11, 2, 1^{2}  $ & $  20558398 $ & & $ 7, 4, 1^{4} $ & $ 5558 $ & &$5, 4^{2}, 1^{2}$ & $362 $ & & $ 4, 2^{3}, 1^{5}$ & $-58$ \\
$ 11, 1^{4}  $ & $ 20024414 $ & & $ 7, 3^{2}, 2 $ & $ 5566 $ & &$5, 4, 3^{2}$ & $350$ & & $4, 2^{2}, 1^{7}$ & $-51$ \\
$ 10, 5   $ & $  -2447421 $ & & $ 7, 3^{2}, 1^{2}$ & $ 5442 $ & &$5, 4, 3, 2, 1$ & $329$& & $4, 2, 1^{9}$ & $-42$ \\
$ 10, 4, 1 $ & $ -2288506 $ & &$ 7, 3, 2^{2}, 1 $ & $ 5246 $ & & $5, 4, 3, 1^{3}$ & $314$& & $4, 1^{11} $ & $-31$ \\
$ 10, 3, 2 $ & $ -2202685 $ & & $ 7, 3, 2, 1^{3} $ & $ 5087 $ & & $5, 4, 2^{3} $ & $302 $& & $ 3^{5} $ & $134$ \\
$ 10, 3, 1^{2}  $ & $   -2169311 $ & & $ 7, 3, 1^{5} $ & $ 4822 $ & &$5, 4, 2^{2}, 1^{2}$ & $294$  & & $3^{4}, 2, 1$  &  $119$  \\
$ 10, 2^{2}, 1  $ & $  -2121105 $ & & $ 7, 2^{4} $ & $ 4854 $ & & $5, 4, 2, 1^{4}$ &  $278$ & & $3^{4}, 1^{3}$ & $110$ \\
$ 10, 2, 1^{3}   $ & $  -2076606 $ & & $ 7, 2^{3}, 1^{2} $ & $ 4766 $ & &$5, 4, 1^{6} $ & $254$ & & $3^{3}, 2^{3}$ & $98$  \\
$ 10, 1^{5} $ & $ -2002441 $ & & $ 7, 2^{2}, 1^{4} $ & $ 4590 $  & &$5, 3^{3}, 1 $ & $290 $ & & $3^{3}, 2^{2}, 1^{2}$ & $94$ \\
$  9, 6  $ & $ 333674 $ & & $ 7, 2, 1^{6} $ & $ 4326 $ & & $5, 3^{2}, 2^{2}$ & $274$& &  $3^{3}, 2, 1^{4}$ & $86$\\
$ 9, 5, 1   $ & $ 293702 $ & & $ 7, 1^{8} $ & $ 3974 $ & & $5, 3^{2}, 2, 1^{2}$ & $266$& &  $3^{3}, 1^{6}$ & $74$\\
$ 9, 4, 2  $ & $ 271934 $ & & $ 6^{2}, 3 $ & $ -3430 $ & & $5, 3^{2}, 1^{4}$ & $250$& &  $3^{2}, 2^{4}, 1$ & $62$ \\
$ 9, 4, 1^{2}  $ & $  266990 $ & & $ 6^{2}, 2, 1 $ & $-3205$ & & $5, 3, 2^{3}, 1$ & $239$& &  $3^{2}, 2^{3}, 1^{3}$ & $59$ \\
$ 9, 3^{2}  $ & $  262226 $ & & $ 6^{2}, 1^{3} $ &  $-3046$  & &$5, 3, 2^{2}, 1^{3}$ & $230$ & &  $3^{2}, 2^{2}, 1^{5}$ & $54$\\
$ 9, 3, 2, 1  $ & $ 254279  $ & & $ 6, 5, 4 $ &  $-1789$  & & $5, 3, 2, 1^{5}$ & $215$& & $3^{2}, 2, 1^{7}$ & $47$ \\
$  9, 3, 1^{3} $ & $  247922 $ & & $ 6, 5, 3, 1 $ & $-1617$  & &$5, 3, 1^{7}$ & $194$ & & $3^{2}, 1^{9}$ & $38$ \\
$   9, 2^{3} $ & $ 244742 $ & & $6, 5, 2^{2} $ &   $-1543$  & &$5, 2^{5}$ & $194$  & & $3, 2^{6}$ & $14$ \\
$   9, 2^{2}, 1^{2} $ & $ 241034 $ & & $6, 5, 2, 1^{2}$ & $ -1501$ & &$5, 2^{4}, 1^{2}$ & $190$ & & $3, 2^{5}, 1^{2} $ & $14$  \\
$  9, 2, 1^{4} $ & $ 233618 $ & & $ 6, 5, 1^{4} $ & $-1417$  & &$5, 2^{3}, 1^{4}$ & $182$& & $3, 2^{4}, 1^{4} $ & $14$ \\
$  9, 1^{6}  $ & $ 222494 $ & & $6, 4^{2}, 1 $ & $-1411$  & &$5, 2^{2}, 1^{6} $ & $170$ & & $3, 2^{3}, 1^{6}$ & $14$ \\
$ 8, 7  $ & $ -65821 $ & & $6, 4, 3, 2 $ & $-1282$  & &$5, 2, 1^{8}$ & $154$& & $3, 2^{2}, 1^{8}$ & $14$ \\
$ 8, 6, 1  $ & $ -49546 $ & &$6, 4, 3, 1^{2} $ & $ -1246$   & &$5, 1^{10} $ & $134$ & & $3, 2, 1^{10}$ & $14$\\
$ 8,5,2 $    &   $-41701$     & & $6,4,2^{2},1$ &  $-1181$      & & $4^{3}, 3$         & $-331$   & & $3, 1^{12}$ & $14$ \\
$ 8, 5, 1^{2}  $ & $  -40775 $ & & $6, 4, 2, 1^{3} $ & $-1141$  & &$4^{3}, 2, 1$ & $-298$ & & $2^{7}, 1$ & $-49$\\
$  8, 4, 3  $ & $ -38146 $ & & $6, 4, 1^{5} $ & $-1066$  & &$4^{3}, 1, 1, 1 $ & $-277$ & & $2^{6}, 1^{3}$ & $-46$ \\
$ 8, 4, 2, 1  $ & $ -36715 $ & & $6, 3^{3} $ & $-1105$  & &$4^{2}, 3^{2}, 1 $ & $-226$ & & $2^{5}, 1^{5}$ & $-41$\\
$ 8, 4, 1^{3}  $ & $  -35602 $ & & $6, 3^{2}, 2, 1$ & $-1054$  & &$4^{2}, 3, 2^{2}$ & $-210$ & & $2^{4}, 1^{7} $ & $-34$  \\
$  8, 3^{2}, 1 $ & $ -34961 $ & & $6, 3^{2}, 1^{3}$ & $-1015$   & & $4^{2}, 3, 2, 1^{2}$ & $-202$& & $2^{3}, 1^{9} $ & $-25$\\
$ 8, 3, 2^{2}   $ & $ -33991 $ & & $6, 3, 2^{3}$ & $ -991$   & & $4^{2}, 3, 1^{4}$ & $-186$& & $2^{2}, 1^{11} $ & $-14$\\
$  8, 3, 2, 1^{2}  $ & $  -33373 $ & & $6, 3, 2^{2}, 1^{2} $ & $-969$   & &$4^{2}, 2^{3}, 1 $ & $-175$ & & $2, 1^{13} $ & $-1$ \\
$   8, 3, 1^{4} $ & $ -32137 $ & & $6, 3, 2, 1^{4}$ & $ -925$    &&$4^{2}, 2^{2}, 1^{3} $ & $-166$ & & $1^{15}$ & $14$\\
$ 8, 2^{3}, 1  $ & $ -31786 $ & & $ 6, 3, 1^{6} $ & $-859$    & &$4^{2}, 2, 1^{5}$ & $-151$ & & $ $   & $ $\\
$ 8, 2^{2}, 1^{3}  $ & $  -30991 $ & & $ 6, 2^{4}, 1 $ & $ -877$  & & $4^{2}, 1^{7}$ & $-130$& & &\\
$  8, 2, 1^{5} $ & $ -29666 $ & & $ 6, 2^{3}, 1^{3} $ & $-850$  & & $4, 3^{3}, 2$ & $-81$& &  &\\
$  8, 1^{7}  $ & $  -27811 $ & & $6, 2^{2}, 1^{5} $ & $ -805$  & & $4, 3^{3}, 1^{2}$ & $-79$& &  &\\
\hline
\end{tabular}
\end{minipage}

}

\end{section}

\end{document}